\documentstyle[amssymb,amsmath,eucal]{article}

\begin{document}

%\righthyphenmin=2

\def\K{{\Bbb K}}
\def\R {{\Bbb R }}
\def\C {{\Bbb C }}
\def\H {{\Bbb H }}
\def\CC {{\Bbb C}^\circ}
\def\RR {{\Bbb R}^\circ}
\def\HH {{\Bbb H}^\circ}

\def\U{{\rm U}}
\def\Sp{{\rm Sp}}
\def\const{{\rm const}}
\def\B{{\rm B}}
\def\O{{\rm O}}
\def\SO{{\rm SO}}
\def\SOS{{\rm SO}^*}
\def\GL{{\rm GL}}

\def\i{{\rm\bf i}}
\def\j{{\rm\bf j}}
\def\k{{\rm\bf k}}

\def\phi{\varphi}
\def\psi{\varpsi}
\def\epsilon{\varepsilon}
\def\le{\leqslant}
\def\ge{\geqslant}

\def\MR{{\cal M}(\R)}
\def\MC{{\cal M}(\C)}
\def\half {{\textstyle\frac12}}
\newcommand{\im}{\mathop{\rm Im}\nolimits}
\newcommand{\re}{\mathop{\rm Re}\nolimits}

\def\kvadrat{\hfill{$\blacksquare$}}

\newcommand{\DET}{\mathop{\rm DET}\nolimits}

\newcommand{\matr}[4]{\bigl(\genfrac{#1}{ #3}\genfrac{ #2}{ #4}\bigr) }

\begin{center}
  {\large\bf  Matrix analogs of the integral 
$$\int_0^\infty \frac{x^{\alpha-1}dx}{(1+x)^\rho}=\B(\alpha,\rho-\alpha)$$
and Plancherel  formula for Berezin kernel representations }

  \medskip
 Neretin Yu.A.\footnote{ Supported by grants RFBR 98-01-00303 
and RFBR 96-01-96249}

         neretin@main.mccme.rssi.ru
\end{center}
\medskip
%\vspace{22pt}

\newcounter{sec}
 \renewcommand{\theequation}{\arabic{sec}.\arabic{equation}}
%\newcounter{fact}
%\newcommand{\fact}[1]{\addtocounter{fact}{1} {\sc #1 \arabic{sec}.\arabic{fact}.}}

The purpose of this paper is to obtain Plancherel formula for Berezin
kernel representations.
We obtain an explicit expression for density of the
 Plancherel measure for all
 10 series of classical Lie groups for sufficiently large values
of parameter of representation
(precisely for the case when representation is square integrable).
A statement of the problem is contained in our Subsection 1.4
(see also recent works on this subject \cite{NO}, \cite{vD},
\cite{OO}, \cite{OZ}, \cite{Ner3},
\cite{Ner4}, \cite{UU}).
 A solution is contained in \S 9.

An evaluation of the Plancherel measure is reduced to an evaluation 
of some integrals over symmetric spaces. 
The main matter of the paper is this evaluation.
In fact we evaluate integrals, which are more general than it is necessary
for the Plancherel formula. Nevertheless neither calculations, nor final formulas
doesn't become more complicated

We give explicit Plancherel formula only for the groups 
  $\O(p,q)$.
For all other series the Plancherel formula can be easily
 obtained in the same way by the a combination of
integrals  (0.1)--(0.11), (0.14)--(0.15), (0.17)--(0.18) 
with classical Gindikin -- Karpelevich formula.

I thanks V.F.Molchanov for discussion of this subject.

\smallskip

{\bf The structure of the paper} is following. Section 0 contains a 
list of integrals that we have evaluated. Evaluations themself are contained
in \S\S4-8. The simplest cases are discussed in \S4.
The  basic case is contained in Subsection 5.1.
The main obstacle is surmounted in Subsection 6.1. All other
 Subsections of \S5--8 contain some technical (but not obvious)
 variations of calculations of 5.1 and 6.1.

In \S\S 1--3 we discuss different models of Riemann noncompact symmetric spaces
(matrix balls, matrics cones, matrix wedges, and sections of wedges).
This information have no relation to calculations themself. We only explain,
what domains of integration in \S\S 4--8 are symmetric spaces.

In \S9 we obtaine the Plancherel formula itself (formula 9.2). 

\medskip

%\vspace{22pt}

{\large \bf \S 0. Formula for Integrals}

%\vspace{22pt}
\medskip

{\bf 0.0. Notations.} By $\K$ we denote  $\R$, $\C$, or the algebra of quaternions 
 $\H$. Denote the quaternionic imaginary units by $\i$, $\j$, $\k$. 
We denote by $\phantom{a}^t$  the operation of transposition of a matrix.
 Denote by $A^*$ the matrix adjoint to a matrix  $A$
 (i.e. $A$ is the matrix obtained from $A$
by conjugation and (complex or quaternionic) conjugation 
of matrix elements). For real matrices  $A^*:=A^t$.
 The symbol  $z^{[g]}$ is introduced by  formula (1.1).

By $[g]_p$ we denote $p\times p$ {\it left upper block
 of a matrix}
$g$.

The symbol  $\det(g)$ denote determinant. Let $g$
be an operator in $\H^n$. Then we  consider $g$ as 
an operator $g_\R:\R^{4n}\to \R^{4n}$
and assume $\det(g)=\sqrt[4]{\det g_\R}$.

The symbol $T\gg0$ means that 
a matrix $T$ is hermitian (i.e. $T=T^*$) positive definite matrix
(i.e. 
 $q T q^*>0$ for any row $q$).

A matrix $R$ is named  {\it dissipative}
if $R+R^*\gg0$. For any dissipative matrix $R$ 
and for any complex  $\tau$ we define quantity
$\det R^\tau=\det(R^\tau)$.
Let $\lambda_1,\dots,\lambda_n$ be eigenvalues of the matrix
  $R$. Then real parts of the complex numbers $\lambda_j$ are positive.
We assume
$$\det R^\tau:=\exp\{\tau\sum\ln\lambda_j\}$$
In the quaternionic case, we define powers of determinant by the formula
$$\det R^\tau=\det R_\R^{\tau/4}.$$

We also use nonstandard notation
$$a^{\{\sigma\|\tau\}}:=a^\sigma\overline a^\tau,
\qquad a,\lambda,\tau\in\C$$
and the similar notation for {\it complex dissipative} matrices $R$:
$$\det R^{\{\sigma\|\tau\}}:=\det R^\sigma\overline {\det R}^\tau$$

The symbol $dT$, where $T$ is matrix, denote the Lebesgue measure on 
 space of matrices. We also fix  following normalization of
Lesbegue measure.
 Let for instance   $T$ pass 
the space of $n\times n$ hermitian matrices. Denote by $t_{kk}$ 
diagonal matrix elements of the matrix $T$, by $t_{km}=p_{km}+iq_{km}$ 
we denote non-diagonal elements. Then we assume
$$dT=
\prod_{k=1}^n dt_{kk}\prod_{1\le k< m \le n} dp_{km}dq_{km}
$$
In the same way, we define a normalization for Lebesgue measure
for spaces of matrices of other types
 (symmetric, antisymmetric etc.)

We start the list of integrals over symmetric spaces evaluated in the paper
 (\S\S 4--8).
Discussion of models of symmetric spaces is contained below in \S\S 1--3.

In all cases $\lambda_j, \sigma_j,\tau_j\in\C$, we also assume
$$
\lambda_{n+1}=0;\qquad \sigma_{n+1}=0;\qquad\tau_{n+1}=0$$
We don't formulate explicitely conditions of convergence since they 
can be easyly observed from explicit
formula.

%\medskip

{\bf 0.1. Integrals over spaces $\GL(n,\R)/\O(n)$, $\GL(n,\C)/\U(n)$,
$\GL(n,\H)/\Sp(n)$.}  Denote by $C_n(\K)$ the space of 
positive defined  $n\times n$ hermitian matrices over  $\K$.
In  \S 4  we establish the following formula:
\begin{eqnarray}
\int_{T\in C_n(\R)}\prod_{j=1}^n
\frac{\det [T]_j^{\lambda_j-\lambda_{j+1}}}
{\det(1+[T]_j)^{\sigma_j-\sigma_{j+1}}}
\det T^{-(n+1)/2}dT= \qquad\qquad\qquad\\
=\prod\limits_{k=1}^n \pi^{(k-1)/2}
\frac{\Gamma(\lambda_k-(k-1)/2)\Gamma(\sigma_k-\lambda_k -(n-k)/2)}
     {\Gamma(\sigma_k-(n-k)/2)}        \nonumber\\
\int_{T\in C_n(\C)}\prod_{j=1}^n
\frac{\det [T]_j^{\lambda_j-\lambda_{j+1}}}
{\det(1+[T]_j)^{\sigma_j-\sigma_{j+1}}}
\det T^{-n}dT=\qquad\qquad\qquad\qquad\\
=\prod\limits_{k=1}^n \pi^{k-1}
\frac{\Gamma(\lambda_k-k+1)\Gamma(\sigma_k-\lambda_k - n+k)}
     {\Gamma(\sigma_k-n+k)}        \nonumber \\
\int_{T\in C_n(\H)}\prod_{j=1}^n
\frac{\det [T]_j^{\lambda_j-\lambda_{j+1}}}
{\det(1+[T]_j)^{\sigma_j-\sigma_{j+1}}}
\det T^{-(2n-1)}dT=\qquad\qquad\qquad\\
=\prod\limits_{k=1}^n \pi^{2k-2}
\frac{\Gamma(\lambda_k-2k+2)\Gamma(\sigma_k-\lambda_k - 2n+2k)}
     {\Gamma(\sigma_k-2n+2k)}        \nonumber
\end{eqnarray}

%\medskip

%\pagebreak

{\bf 0.2.
 Integrals over the spaces
 $\O(n,n)/\O(n)\times\O(n)$, $\U(n,n)/\U(n)\times\U(n)$,
$\Sp(n,n)/\Sp(n)\times\Sp(n)$.}
 We denote by $W_n(\K)$ the space of $n\times n$
dissipative matrices $R$ over $\K$. We represent a dissipative matrix  $R$
 in the form $R=T+S$ where  $T=T^*\gg0$,  $S^*=-S$.
In  \S 5  we evaluate the following integrals:
\begin{eqnarray}
\int_{T+S\,\in W_n(\R)}
\prod_{j=1}^n
\frac{ \det[T]_j^{\lambda_j-\lambda_{j+1}}    }
     { \det[1+T+S]_j^{\sigma_j-\sigma_{j+1}}  }
    \det T^{-n} dT\,dS= \qquad\qquad\\
= \prod_{k=1}^n
\frac{  \pi^{k-1}
 \Gamma(\lambda_k-(n+k)/2+1)\Gamma(\sigma_k-\lambda_k - (n-k)/2)  }
 {  \Gamma(\sigma_k-n+k) } \nonumber\\
\int_{T+S\,\in W_n(\C)}
\prod_{j=1}^n
\frac{ \det[T]_j^{\lambda_j-\lambda_{j+1}}    }
     { \det[1+T+S]_j^{\{\sigma_j-\sigma_{j+1}\|\tau_j-\tau_{j+1}\}}  }
    \det T^{-2n} dT\,dS=\\
=\prod_{k=1}^n
2^{2(n-k+1)-\sigma_k-\tau_k}\pi^{2k-1}
 \frac{ \Gamma(\lambda_k-n-k+1)\Gamma(\sigma_k+\tau_k-\lambda_k-n+k)}
  {\Gamma(\sigma_k-n+k)    \Gamma(\tau_k-n+k)    }\nonumber\\
\int_{T+S\,\in W_n(\H)}
\prod_{j=1}^n
\frac{ \det[T]_j^{\lambda_j-\lambda_{j+1}}    }
     { \det[1+T+S]_j^{\sigma_j-\sigma_{j+1}}  }
    \det T^{-4n} dT\,dS=\qquad\qquad\\
= \prod_{k=1}^n
2^{\sigma_k-4(n-k+1)}\pi^{4k-3}
\frac{  \Gamma(\lambda_k-2(n+k)+1)\Gamma(\sigma_k-\lambda_k-2(n-k))  }
    { \Gamma(\sigma_n/2-2(n-k))\Gamma (\sigma_n/2-2(n-k)-1) } \nonumber
\end{eqnarray}

\medskip

{\bf 0.3.
 Integrals over the spaces
 $\O(p ,q)/\O(p)\times\O(q)$, $\U(p,q)/\U(p)\times\U(q)$,
$\Sp(p,q)/\Sp(p)\times\Sp(q)$.}
Let $p\le q$. Let us define a space $SW_{p,q}(\K)$. Its points are pairs
 $(K,L)$, where $K$ is a $p\times p$ matrix,  $L$ is
a $p \times (q-p)$ matrix, and
$$ \frac12 (K+K^*)-LL^*\gg0$$
We represent $K$ in the form $K=M+N$, where $M=M^*$, $N=~-N^*$.
In \S 6  we evaluate the following integrals
\begin{eqnarray}
\int_{SW_{p,q}(\R)}
\prod_{j=1}^p
\frac{ \det[M-LL^*]_j^{\lambda_j-\lambda_{j+1}}  }
  {\det [1+M+N]_j^{\sigma_j-\sigma_{j+1}} }
  \cdot
\det(M-LL^*)^{-(p+q)/2} dM\,dN\,dL =\\
=\prod_{k=1}^p
\pi^{k-(q-p)/2-1}
\frac{\Gamma(\lambda_k-(q+k)/2+1)\Gamma(\sigma_k-\lambda_k-(p-k)/2)}
     {\Gamma(\sigma_k-p+k)}\nonumber\\
\int_{SW_{p,q}(\C)}
\prod_{j=1}^p
\frac{ \det[M-LL^*]_j^{\lambda_j-\lambda_{j+1}}  }
  {\det [1+M+N]_j^{\{\sigma_j-\sigma_{j+1}\|\tau_j-\tau_{j+1}\}} }
  \times\qquad\qquad\qquad\qquad\qquad\quad \nonumber \\  \times
\det(M-LL^*)^{-(p+q)} dM\,dN\,dL =\\
=           \prod_{k=1}^p
2^{2(p-k)-\sigma_k-\tau_k}\pi^{2k-1}
\frac{ \Gamma(\lambda_k-(q+k)+1)\Gamma(\sigma_k+\tau_k-\lambda_k-p+k)}
     {\Gamma(\sigma_k-p+k)  \Gamma(\tau_k-p+k)}\nonumber\\
\int_{SW_{p,q}(\H)}
\prod_{j=1}^p
\frac{ \det[M-LL^*]_j^{\lambda_j-\lambda_{j+1}}  }
  {\det [1+M+N]_j^{\sigma_j-\sigma_{j+1}} }
  \cdot
\det(M-LL^*)^{-2(p+q)} dM\,dN\,dL =\\
=\prod_{k=1}^p 2^{4(p-k+1)-\sigma_k}\pi^{2(q+k)-3}
  \frac{\Gamma(\lambda_k-2(q+k)+1)\Gamma(\sigma_k-\lambda_k-2(p-k))}
     {\Gamma(\sigma_k/2-2(p-k))    \Gamma(\sigma_k/2-2(p-k)-1)}
\nonumber
\end{eqnarray}
(we assume $\lambda_{p+1}=\sigma_{p+1}=\tau_{p+1}=0$).

\medskip

{\bf 0.4. Integrals over the spaces $\Sp(2n,\R)/\U(n)$
and $\Sp(2n,\C)/\Sp(n)$.}
We realize the symmetric space   $\Sp(2n,\R)/\U(n)$
as the space of all symmetric   $n\times n$ matrices $R$ 
 with a positive definite real part. Let us represent $R$  in the form $R=T+iS$, 
where
matrices $T,S$ are real.   In \S7 we evaluate the following integral
\begin{eqnarray}
\int\limits_{T=T^t\gg0,\,S=S^t}
\prod_{j=1}^n
\frac{ \det[T]_j^{\lambda_j-\lambda_{j+1}}    }
     { \det[1+T+iS]_j^{\{\sigma_j-\sigma_{j+1}\|\tau_j-\tau_{j+1}\}}  }
    \det T^{-(n+1)} dT\,dS=  \\
= \prod_{k=1}^n
\frac{2^{2-\sigma_k-\tau_k+n-k}  \pi^{k}
 \Gamma(\lambda_k-(n+k)/2)\Gamma(\sigma_k+\tau_k-\lambda_k - (n-k)/2)  }
 {  \Gamma(\sigma_k-(n-k)/2) \Gamma(\tau_k-(n-k)/2)} \nonumber
\end{eqnarray}

We realize the space $\Sp(2n,\C)/\Sp(n)$ as the space of pairs $(T,S)$
of complex  $n\times n$ matrices satisfying the conditions
$$T=T^*\gg0,\qquad S=S^t$$
Let $\j$ be the quaternionic imaginary unit. Then
\begin{eqnarray}
\int\limits_{T=T^*\gg0,\,S=S^t}
\prod_{s=1}^n
\frac{ \det[T]_s^{\lambda_s-\lambda_{s+1}}    }
     { \det[1+T+ S\j]_s^{\sigma_s-\sigma_{s+1}}  }
    \det T^{-(2n+1)} dT\,dS=\qquad \\
=\prod_{k=1}^n 2\pi^{2k-1}
 \frac{\Gamma(\lambda_k-(n+k))\Gamma (\sigma_k-\lambda_k-(n-k))}
        {\Gamma (\sigma_k-2(n-k)-1)}  \nonumber
\end{eqnarray}

\medskip

{\bf 0.5. Integrals over the spaces $\O(n,\C)/\O(n)$, $\SOS(2n)/\U(n)$.}
Let
\begin{equation}
J=\begin{pmatrix}0& -1& &&       \\
               1&  0& &&         \\
                &&      0& -1 &    \\
                &&      1&  0 &      \\
                &&    &&  & \ddots
\end{pmatrix}\qquad
\end{equation}
We realize the space $\O(2n,\C)/\O(2n)$ (respectively $\SOS(4n)/\U(2n)$)
 as the space of  
 $2n\times 2n$ matrices $R$ over $\R$
 (respectively over $\C$),
satisfying the conditions
\begin{equation}
J^{-1}R^tJ= R;\qquad R+R^*\gg0
\end{equation}
Assume $R=T+S$, where $T=T^*,S=-S^*$. Then,
in the case    $\O(2n,\C)/\O(2n)$,    we have
\begin{eqnarray}
\int
\prod_{j=1}^n
\frac{ \det[T]_{2j}^{(\lambda_j-\lambda_{j+1})/2}    }
     { \det[1+T+S]_{2j}^{(\sigma_j-\sigma_{j+1})/2}  }
    \det T^{-(2n-1)/2} dT\,dS=\qquad\qquad \\ \qquad\qquad=
 \prod_{k=1}^n \pi^{2k-2}
  \frac{\Gamma(\lambda_k-(n+k)+2)\Gamma(\sigma_k-\lambda_k-(n-k))}
   {\Gamma(\sigma_k-(n-k))}   \nonumber
\end{eqnarray}
In the case $\SOS(4n)/\U(2n)$,    we obtain
\begin{eqnarray}
\int
\prod_{j=1}^n
\frac{ \det[T]_{2j}^{(\lambda_j-\lambda_{j+1})/2}    }
     { \det[1+T+S]_{2j}^{\{(\sigma_j-\sigma_{j+1})/2\|(\tau_j-\tau_{j+1})/2\}}  }
    \det T^{-(2n-1)} dT\,dS=    \qquad\qquad       \\=
 \prod_{k=1}^n
 2^{2+4(n-k)-\sigma_k-\tau_k}\pi^{4k-3}
\frac{\Gamma(\lambda_k-2(n+k)+3)\Gamma(\sigma_k+\tau_k-\lambda_k-2(n-k))}
         {\Gamma(\sigma_k-2(n-k))\Gamma(\tau_k-2(n-k))}\nonumber
\end{eqnarray}

We realize spaces $\O(2n+1,\C)/\O(2n+1)$ and $\SOS(4n+2)/\U(2n+1)$
as spaces of 
$(1+1+2n) \times (1+1+2n) $ matrices $R$, having block structure
\begin{equation}    R=
\begin{pmatrix}
1&0&*\\ 0&1&0\\0&*&*
\end{pmatrix},
\end{equation}
and satisfying condition (0.13). Then
in the case   $\O(2n+1,\C)/\O(2n+1)$    we have
\begin{eqnarray}
\int
\prod_{j=1}^n
\frac{ \det[R+R^t]_{2j+2}^{(\lambda_j-\lambda_{j+1})/2}    }
     { \det[1+R]_{2j+2}^{(\sigma_j-\sigma_{j+1})/2}  }
    \det (R+R^*)^{-n} dR=\qquad\qquad \\  \qquad\qquad =
C_1 \prod_{k=1}^n \pi^{2k-1}
  \frac{\Gamma(\lambda_k-(n+k)+1)\Gamma(\sigma_k-\lambda_k-(n-k))}
   {\Gamma(\sigma_k-(n-k))}   \nonumber
\end{eqnarray}
and in the case $\SOS(4n+2)/\U(2n+1)$  we have
\begin{eqnarray}
\int
\prod_{j=1}^n
\frac{ \det[R+R^*]_{2j+2}^{(\lambda_j-\lambda_{j+1})/2}    }
     { \det[1+R]_{2j+2}^{\{(\sigma_j-\sigma_{j+1})/2\|(\tau_j-\tau_{j+1})/2\}}  }
    \det (R+R^*)^{-2n} dR= \qquad\qquad          \\ \qquad\qquad=
C_2 \prod_{k=1}^n
\pi^{4k-1}
\frac{\Gamma(\lambda_k-2(n+k)+2)\Gamma(\sigma_k+\tau_k-\lambda_k-2(n-k))}
         {\Gamma(\sigma_k-2(n-k))\Gamma(\tau_k-2(n-k))}\nonumber
\end{eqnarray}
where $C_1$, $ C_2$ are some powers of 2.

\medskip

{\bf 0.5. Known partial cases.}
a)  Hua Loo Keng integrals (\cite{Hua}) are integrals over matrix balls
(see below \S 1) having type
$$\int\det(1-zz^*)^\alpha dz$$
They can be reduced to very partial cases of our integrals by means of the Cayley transform $z=(1+R)^{-1}(1-R)$.

b) Integrals (0.1)--(0.3) for $\lambda_1=\dots=\lambda_n$ can be reduced
by substitution $T+1=(H-1)^{-1}$ to a partial case of 
Gindikin  $\B$-function of cone,
see \cite{Gin}, see also \cite{FK},
Theorem VII.1.7. 

c) Integrals (0.1)--(0.3) for $\sigma_1=\dots=\sigma_n$ were evaluated
by Faraut and Koranyi (see. \cite{FK}, XIV.4).

d)Berezin in one of the last his papers (\cite{Ber2})
  announced the Plancherel formula for kernel representations
of the groups 
 $\U(p,q)$, $\Sp(2n,\R)$, $\SOS(2n)$.
In fact, it is equivalent to an evalution of integrals  (0.5), (0.8),
 (0.10), (0.15), (0.18)
for $\sigma_1=\tau_1=\sigma_2=\tau_2=\dots$
(up to meromorphic multiplier
depending on $\sigma$)%
\footnote{after submitting of the paper I received interesting preprint of G.Zhang \cite{Zha}
generalizing \cite{UU}}. Proofs were published by
Upmeier and Unterberger \cite{UU}. Our method  differs
from 
\cite{UU} and (as far as I know)
it differs from Berezin method.

e) Plancherel formula is known for rank 1 classical groups, i.e. $\O(n,1)$, $\Sp(n,1)$,
see. \cite{vD},
This case is essentially more simple than the case of rank $\ge2$.
(our calculations in this case are trivial).

\medskip
\vspace{22pt}

{\large\bf \S 1. Models of Symmetric Spaces: Matrix Balls}

\nopagebreak
\vspace{22pt}
\nopagebreak

\addtocounter{sec}{1}
\setcounter{equation}{0}

\medskip

{\bf 1.1. Matrix balls.}
By a {\it norm}   of a matrix
we always mean the norm of a operator from a coordinate
Euclidean space 
$\K^p$ to $\K^q$.
By a {\it matrix ball} we mean

\smallskip

the set of all
$p\times p$ matrices over $\K$ having a norm $<1$ (where $\K=\R,\C,\H$)

\smallskip

or set of all $n\times n$ matrices over $\K$,
satisfying one of symmetry conditions: {\it
symmetric} ($z=z^t$),
{\it antisymmetric} ($z=-z^t$), {\it hermitian} ($z=z^*$),
{\it anti-hermitian} ($z=-z^*$).

\smallskip

Let  $G$ be one of classical groups $\GL(n,\R)$, $\GL(n,\C)$,
$\GL(n,\H)$, $\O(p,q)$,  $\U(p,q)$, $\Sp(p,q)$, $\Sp(2n,\R)$,
$\Sp(2n,\C)$,  $\O(n,\C)$, $\SOS(2n)$. Let $K$ be
a maximal compact subgroup in $G$.
Any  non-compact classical Riemann symmetric space
$G/K$ has  matrix ball type
realization (see \cite{Ner1}).
 Let us enumerate these realizations. Emphasis
that {\it in all cases we consider a space of matrices with norm} $<1$.

\smallskip 

 1. $\O(p,q)/\O(p)\times\O(q)$ is
the space  of
$p\times q$ matrices over $\R$.

\smallskip 

 2. $\GL(n,\R)/\O(n)$ is the space
of symmetric matrices   $n\times n$ matrices over $\R$.

\smallskip

 3. $\O(n,\C)/\O(n)$ -- is the space
 of antisymmetric  
 $n\times n$ matrices over $\R$.

\smallskip

 4. $\U(p,q)/\U(p)\times\U(q)$ 
is the space of  
$p\times q$ matrices over $\C$.

\smallskip

 5. $\Sp(2n,\R)/\U(n)$ is the space of symmetric 
 $n\times n$ matrices over $\C$

\smallskip

 6. $\SO^*(2n)/\U(n)$ is the space of antisymmetric  
 $n\times n$ matrices over $\C$

\smallskip

 7. $\GL(n,\C)/\U(n)$ is the space of hermitian  
 $n\times n$ matrices over $\C$

\smallskip

 8. $\Sp(p,q)/\Sp(p)\times\Sp(q)$ is the space 
 $p\times q$ matrices over $\H$.

\smallskip

 9. $\GL(n,\H)/\Sp(n)$    is the space of hermitian $n\times n$
matrices  $\H$

\smallskip

10. $\Sp(2n,\C)/\Sp(n)$ is the space of $n\times n$ anti-hermitian matrices
over $\H$ 

\smallskip

In all cases, the group  $G$ acts on the matrix ball
by fractional linear transformations having the form
\begin {equation}
g:z\mapsto z^{[g]}:=(\alpha+z\gamma)^{-1}(\beta+z \delta)
\end{equation}
Moreover, the group  $G$ consist of all fractional linear transformations,
mapping bijectively matrix ball to itself. Let us describe the group $G$
explicitly.

In the cases  1,4,8 (i.e. in the case $G=\U(p,q,\K)$) we
consider the space 
$\K^p\oplus\K^q$ equipped with indefinite hermitian form   $M=
\left(\begin{array}{cc} 1&0\\ 0&-1 \end{array} \right)  $.
We realize the pseudo-unitary group  $\U(p,q,\K)$ over $\K$ (i.e. one of the groups
$\O(p,q)$, $\U(p,q)$, $\Sp(p,q)$) as the group of block matrices
$g=\left(\begin{array}{cc}\alpha &\beta\\ \gamma&\delta \end{array} \right)$,
having the size $(p+q)\times(p+q)$ and preserving the hermitian form $M$.

In all other cases, we consider 
the space
$\K^n\oplus\K^n$ equipped with indefinite hermitian form $M=
\left(\begin{array}{cc} 1&0\\ 0&-1 \end{array} \right)$.
As above elements of the group $G$ preserve this form
and also they preserve once more form $\Lambda$
depending on type of the matrix space:

\smallskip

a) in the case of a space of symmetric matrices (
$G/K=$
$\GL(n,\R)/\O(n)$,
$\Sp(2n,\R)/\U(n)$), the form $\Lambda$ is an antisymmetric 
bilinear form with the matrix
$\left(\begin{array}{cc} 0&1\\ -1&0 \end{array} \right)$.

\smallskip

b) in the case of a space of antisymmetric matrices
($G/H=\O(n,\C)/\O(n,\R)$,  $\SO^*(2n)/\U(n)$),
the form $\Lambda$ is a symmetric  bilinear form with the  matrix
$\left(\begin{array}{cc} 0&1\\ 1&0 \end{array} \right)$.

\smallskip

c) in the case of a space of hermitian matrices
( $G/K=$
 $\GL(n,\R)/\O(n)$, $\GL(n,\C)/\U(n)$, $\GL(n,\H)/\Sp(n)$),
the form $\Lambda$ is an anti-hermitian form with  the matrix
$\left(\begin{array}{cc} 0&1\\- 1&0 \end{array} \right)$;
remind that a sesquilinear form  $\Lambda$ is named {\it anti-hermitian}
if $\Lambda(v,w)=-\overline{\Lambda(w,v)}$ for all $v$, $w$.

\smallskip

d) in the case of a space of anti-hermitian matrices ( $G/H=\Sp(2n,\C)/\Sp(n)$),
the form $\Lambda$ is a hermitian form with the matrix
$\left(\begin{array}{cc} 0&1\\1&0 \end{array} \right)$.

\smallskip

In all cases,
a stabilizer of  the point $z=0$
consists of matrices 
$\left(\begin{array}{cc}\alpha &0\\ 0&\delta \end{array} \right)\in G$,
where matrices $\alpha$  and $\delta$ are unitary over $\K$; note also that
in all cases,
exept $G=\U(p,q,\K)/\U(p,\K)\times\U(q,K)$, we have $\delta=\alpha^{*-1}   $.

\medskip

{\bf 1.2. Jacobian.}  {\sc  Lemma 1.1.} {\it A differential of a map
$z\mapsto z^{[g]}$ in a point   $z$ is }
$$(\alpha+z\gamma)^{-1}dz(-\gamma z^{[g]}+\delta)$$

{\sc Proof.} We have to calculate the following expression up to 
 $o(\epsilon)$  
\begin{multline*}
(\alpha+(z+\epsilon u)\gamma)^{-1}(b+(z+\epsilon u)\delta)=\\
=(\alpha+z\gamma)^{-1}
(1+\epsilon u \gamma (\alpha+z\gamma)^{-1})^{-1}
(\beta+z\delta+\epsilon u\delta)=\\
=(\alpha+z\gamma)^{-1}
(1-\epsilon u \gamma (\alpha+z\gamma)^{-1}+o(\epsilon))
(\beta+z\delta+\epsilon u\delta)=\\
=(\alpha+z\gamma)^{-1}(\beta+z\delta)+
\epsilon (\alpha+z\gamma)^{-1}u
\left[-\gamma (\alpha+z\gamma)^{-1}
(\beta+z\delta) + \delta\right]+o(\epsilon) =\\
=z^{[g]}+\epsilon (\alpha+z\gamma)^{-1}
 u (-\gamma z^{[g]}+  \delta)+o(\epsilon)\qquad \blacksquare
\end{multline*}

Recall a formula for the determinant of block matrix  (see \cite{Gan},\S II.5).
 Below it is used very intensively
\begin{equation}
\det
\left(\begin{array}{cc} A&B\\C&D \end{array} \right)=
\det A\cdot \det (D-CA^{-1}B)
\end{equation}

{\sc  Lemma 1.2.}
 \begin{equation}
\det(-\gamma z^{[g]}+\delta)=(\det g)\det(\alpha+z \gamma )^{-1}
\end{equation}

{\sc Proof.}
$$
\det(-\gamma z^{[g]}+\delta)=
\det(-\gamma (\alpha+z\gamma)^{-1}(\beta+z\delta)+\delta)
$$
By formula (1.2), we reduce this expression to the
form 

\begin{multline*}\det(\alpha+z\gamma)^{-1}
\det
\left(\begin{array}{cc} \alpha+z\gamma&\beta +z\delta\\ \gamma&\delta \end{array} \right)=
\\
=\det(\alpha+z\gamma)^{-1}  \det
\left(\begin{array}{cc} \alpha&\beta \\ \gamma&\delta \end{array} \right)
\qquad\blacksquare
\end{multline*}

Lemma 1.1, 1.2 impliy the following statement

{\sc  Lemma 1.3.} {\it A Jacobian of a transformation
$g\mapsto z^{[g]}$ in a point $z$ is
$$|\det (\alpha+z\gamma) |^{-h}$$
where $h=(p+q)\dim \K$ in the cases} 1,4,8,  {\it and $h=2n\dim G/K$
 in other cases. }

\medskip

{\bf 1.3. Invariant measure.}

{\sc  Lemma 1.4.}
\begin{equation}
1-z^{[g]} (z^{[g]})^*=
(\alpha+z\gamma)^{-1}(1-zz^*)(\alpha+z\gamma)^{*-1}
\end{equation}

{\sc Proof.}
\begin{eqnarray}
1-z^{[g]} (z^{[g]})^*=
1-      (\alpha+z\gamma)^{-1}(\beta+z\delta)
(\beta^*+\delta^*z^*)(\alpha^*+\gamma^*z^*)^{-1}=\nonumber\\
=(\alpha+z\gamma)^{-1} \bigl\{
(\alpha+z\gamma) (\alpha^*+\gamma^*z^*) -(\beta+z\delta)(\beta^*+\delta^*z^*)
\bigr\}  (\alpha^*+\gamma^*z^*)^{-1}
\end{eqnarray}
By
$$\left(\begin{array}{cc} \alpha&\beta \\ \gamma&\delta \end{array} \right)
\left(\begin{array}{cc} 1&0 \\0&-1 \end{array} \right)
\left(\begin{array}{cc} \alpha&\beta \\ \gamma&\delta \end{array} \right)^*
                                   =
\left(\begin{array}{cc} 1&0 \\0&-1 \end{array} \right)
 $$
we obtain that the factor in (1.5) in the curly brackets has required form.
\shoveright{$\blacksquare$}

{\sc  Corollary 1.5.}
$$\det(1-z^{[g]} (z^{[g]})^*)=|\det (\alpha+z\gamma) |^{-2}  \det(1-zz^*)$$

{\sc Proposition 1.6.} {\it A $G$-invariant measure on a matrix ball 
 is given by the formula
$$|\det(1-zz^*)|^{-h/2} dz$$
where $h$ is the same as in Lemma {\rm 1.3}, and $dz$ denote
 Lebesgue measure on the space of matrices.}

{\sc Proof.} It is consequence of Lemma 1.3 and Corollary 1.5.

\medskip

{\bf 1.4. Kernel representations.} Consider a function on the matrix
ball 
$G/K$ given by the formula
$${\cal B}_\alpha(z)=\det(1-zz^*)^{-\alpha}$$
We can consider the function $\cal B$ as a function on the 
the group $G$
$${\cal B}_\alpha(g) :={\cal B}_\alpha(0^{[g]})$$
{\it Kernel representation} $T_\alpha$ of the group $G$ is the unitary
representation
containing a cyclic vector $v$ such that 
 $$\langle T_\alpha (g)v,v\rangle= {\cal B}_\alpha(g)$$.

Representation $T_\alpha$ exists iff the function
${\cal B}_\alpha(g)$ on the group $G$ is positive definite.
Conditions of positive definiteness for
 ${\cal B}_\alpha(g)$ were obtained by Berezin
in \cite{Ber1} (some additions were given by Gindikin \cite{Gin2}).
 For instance, Berezin condition for the groups  $\O(p,q)$ 
 is
$$2 \alpha= 0,1,2,\dots, p-1\qquad \mbox{or} \qquad 2\alpha>p-1$$
(we assume $p\le q$).
For other series conditions are similar.

Kernel representation also can be defined as restrictions of highest
weight representation of some 'overgroup' $G^*\supset G$ to $G$.
They were 
introduced by Berezin \cite{Ber2} for hermitian symmetric spaces, 
in \cite{OO}, \cite{Ner3} there were observed
that the problem is general for all Riemannian noncompact symmetric spaces
(besides several exceptional spaces). Detail discussion and more 'material' 
definitions see \cite{Ner3}, \cite{Ner4}, see also recent works
\cite{NO}, \cite{OO}, \cite{UU}, \cite{vD}.

\vspace{22pt}

{\large\bf \S2. Models of Symmetric Spaces:
 Symmetric Cones and Symmetric Wedges}

\vspace{22pt}

\addtocounter{sec}{1}
\setcounter{equation}{0}

Matrix ball models are not convenient for our calculations.
Our  purpose in \S2-\S3 is to obtain realizations at which  parabolic subgroups
acts by affine transformations.
Our realizations are similar to realizations from     \cite{Pya}.

\medskip

{\bf 2.1. Cayley transform.} The Cayley transform is the transform of the space
of square matrices given by the formula

\begin{eqnarray}
z=\frac{1-R}{1+R}=-1+(1+R)^{-1}\\
R=\frac{1-z}{1+z}=-1+(1+z)^{-1}
\end{eqnarray}

By Lemma 1.1, the differential of the Cayley transform is given by the formula
$$dR=2(1+z)^{-1} dz\,(1+z)^{-1}$$

{\sc Proposition 2.1.} {a)} $zz^*=1\quad \Longleftrightarrow\quad R=-R^*$

{b)} $\|z\|< 1 \quad \Longleftrightarrow\quad   R+R^*\gg0$

 c) $z=z^t \quad \Longleftrightarrow\quad   R=R^t$

 d) $z=z^* \quad \Longleftrightarrow\quad R=R^*$

(Statements a),b) are well-known, and c),d) are obvious.)

\medskip

{\bf 2.2. Symmetric cones.}
 Proposition 2.1b),d)  implies

{\sc Corollary 2.2.} {\it The Cayley transform maps the matrix ball
\begin{multline*}\GL(n,\K)/\U(n,\K)=\\=
\GL(n,\R)/\O(n),\quad \GL(n,\C)/\U(n),\quad \GL(n,\H)/\Sp(n)
\end{multline*}
 to the cone of positive definite hermitian matrices over $\K=\R,\C,\H$
respectively.}

Thus we obtaine another realization for symmetric spaces
enumerated in the Corollary.
  The group  $G=\GL(n,\K)$ acts on the cone of positive definite matrices by the formula
$$A:T\mapsto ATA^*.$$

\medskip

{\bf 2.2. Wedge of dissipative matrices.} Proposition 2.1b)
also implies

{\sc Corollary 2.3.}
 {\it The Cayley transform maps matrix ball
\begin{multline*}\U(n,n,\K)/\U(n,\K)\times \U(n,\K) = \\  =
\O(n,n)/\O(n)\times\O(n),\quad \U(n,n)/\U(n)\times \U(n),\quad
\Sp(n,n)/\Sp(n)\times \Sp(n)
\end{multline*}
to the space of all dissipative matrices over $\K=\R,\C,\H$ respectively.}

In this case, it is convenient to realize the group  $G=\U(n,n,\K)$ as
the group of matrices
$\left(\begin{array}{cc} a&b \\ s&d\end{array} \right)$
preserving
the hermitian form
$\left(\begin{array}{cc} 0&1 \\ 1&0\end{array} \right)$.

The group $G$ acts on the space of dissipative matrices
 by fractional linear transformations
\begin{equation}R\mapsto R^{[g]}:=(a+Rc)^{-1}(b+Rd)\end{equation}

\medskip

{\bf 2.4. Siegel upper half-plane.}
Recall that the matrix ball  $\Sp(2n,\R)/\U(n)$
 consists of complex symmetric matrices.
By Proposition  2.1 b),c) its image under the Cayley transform is the space
of complex symmetric matrices with positive definite real part.

\medskip

{\bf 2.5. The space $\Sp(2n,\C)/\Sp(n)$.} In this case the matrix ball
consists of quaternionic anti-hermitian matrices.
Consider the modified Cayley transform
$$R=-1+2(1+\i z)^{-1}$$

The image of this transform consist of quaternionic matrices $R$
satisfying the conditions
\begin{equation}
R^*=\i^{-1}R\i;  \qquad R+R^*\gg0.
\end{equation}
Represent a matrix $R$  in the form $R=S+T\j$, where $R,S$ are complex matrices.
Condition (2.4) is equivalent to
$$ T=T^*;\qquad S=S^t;\qquad T\gg0.$$

\medskip
{\bf 2.6.  Spaces $\O(2n,\C)/\O(2n,\R), \SOS(4n)/\U(2n)$.}
Recall that these matrix balls consist of 
 $2n\times 2n$ antisymmetric matrices over $\R$, $\C$ respectively.
      Define the matrix
\begin{equation}
J=\begin{pmatrix}0& -1& &&       \\
               1&  0& &&         \\
                &&      0& -1 &    \\
                &&      1&  0 &      \\
                &&    &&  & \ddots
\end{pmatrix}\qquad
\end{equation}
Consider the modified Cayley transform 
$$R=-1+(1+JZ)^{-1}$$
The image of the matrix ball consists of matrices  $R$ satisfying the conditions
\begin{equation}
R^t=JRJ^{-1};\qquad R+R^*\gg0
\end{equation}
We will continue a discussion of these spaces in \S8.
\medskip

{\bf 2.7. Invariant measure.}
{\sc Lemma 2.4.} {\it Let $z$ and $R$
 are linked by the Cayley transform {\rm (2.1)}. Then }
$$1-zz^*=2(1+R)^{-1}(R+R^*)(1+R^*)^{-1} $$
\begin{multline*}        \mbox{\sc Proof.}
1-zz^*=1-(1+R)^{-1}(1-R)(1-R^*)(1+R^*)^{-1}=\\
=(1+R)^{-1}\left\{(1+R)(1+R^*)-(1-R)(1-R^*)\right\}(1+R^*)^{-1}
\end{multline*}
and we remove brackets in the middle factor.
\kvadrat

The Lemma implies two corollaries. In the thirst place, it
allows to rewrite functions
${\cal B}_\alpha$ (see Subsection 1.4) in new coordinates. In all cases
 (see Subsections 2.1--2.6) we obtain
\begin{equation}
{\cal B}_\alpha(R)=\biggl(\frac {\det(2(R+R^*))} {|\det(1+R)|^2}\biggr)^{\alpha}
\end{equation}
In the second place, Lemma 2.4 and Proposition 1.6 give
formula for   $G$-invariant measure on a cone or a wedge.
In all  enumerated cases the invariant measure is given by the formula
\begin{equation}\det (R+R^*)^{\frac{\dim (G/K)}{m}} dR
\end{equation}
where  $m$ is the size of matrices, which we consider.
\medskip

{\bf 2.8. Action of a parabolic subgroup.}
In the case of symmetric cones  $\GL(n,\K)/\U(n,\K)$, we consider
the minimal parabolic subgroup ${\cal P}\subset\GL(n,\K)$ consisting of
lower-triangular matrices.  Obviously, functions
\begin{equation}
\prod_{j=1}^n \det[T]_j^{\lambda_j-\lambda_{j+1}}
\end{equation}
are eigenfunctions of the subgroup $P$.

Consider the case when a space $G/K$ is a wedge of dissipative matrices
(i.e  $G/K=\U(n,n,\K)/\U(n,K) \times \U(n,K)$),
Let us define the group $\cal P$ consisting 
of all transformations 
$$R\mapsto QRQ^*+ S$$
of the wedge, where $Q\in\GL(n,\K)$ is a lower-triangular matrix and $S=-S^*$.
It is easy to see, that  $\cal P$ is a minimal parabolic subgroup
in $G=\U(n,n,\K)$.

Eigenfunctions of the parabolic subgroup  $\cal P$ on the wedge have the form
(2.9) where  $T=(R+R^*)/2$.

For all other wedges the picture is similar.

\medskip

{\bf 2.9. Comments to exterior of 
 integrals \S 0.} Now we are ready to explain a sense of some factors
in integrands (0.1)--(0.18).
 The last factor is  density
of the invariant measure,
see (2.8),
it is possible to include this factor to numerator  
of the fraction).
The numerator of  the fraction
$\prod_{j=1}^n \det[T]_j^{\lambda_j-\lambda_{j+1}}$
is an eigenfunction of the minimal parabolic subgroup.

 I don't understand the sense of the product in denominator.
In any case  (see (2.7),(2.8))
 all expressions of the form
${\cal B}_\alpha(R)\cdot\prod_{j=1}^n \det[T]_j^{\lambda_j-\lambda_{j+1}} $
are contained in the set of our integrands.

\vspace{22pt}

{\large\bf \S 3. Models of Symmetric Spaces: Sections of Wedges}

\vspace{22pt}

\addtocounter{sec}{1}
\setcounter{equation}{0}

In this section we discuss symmetric spaces $G/K$ which have no matrix cone and matrix wedge realizations(it is the case when the group $G$ is not split).

\medskip

{\bf 3.1. The case $G=\U(p,q)/\U(p)\times\U(q)$.}
Denote by $\B_{\alpha,\beta}=\B_{\alpha,\beta}(\K)$
 the space of all 
$\alpha\times \beta$ matrices over $\K$  with norm $<1$. Recal that we identifies
$\B_{\alpha,\beta} (\K)$ with symmetric space
$\U(\alpha,\beta,\K)/\U(\alpha,\K)\times\U(\beta,\K)$, see Subsection 1.1.

Assume $p<q$. Let us wright $Z\in\B_{p,q}$  as block 
$p\times ((q-p))+p)$ matrix:
$$ Z=
\left(\begin{array}{cc} X&Y\end{array}\right)
$$
Complete the matrix $Z$ to the matrix $\widetilde Z$ of the size
 $((q-p))+p)\times ((q-p))+p)$ by the formula
\begin{equation}
\widetilde Z =
\left(\begin{array}{cc}0&0\\ X&Y\end{array}\right)
\end{equation}
The Cayley transform (2.1) maps this matrix $\widetilde Z$
to a dissipative matrix   $R$ having the form
\begin{equation}                  R=
\left(\begin{array}{cc}1&0\\2L&K\end{array}\right)
\end{equation}
where
$$
K=-1+2(1+Y)^{-1};\qquad
L=-(1+Y)^{-1}  X
$$
Conversely, it is easy to check, that 
any dissipative matrix having form (3.2) is a Cayley transform
of some matrix having form (3.1).

Thus we realized the space
$\U(p,q,\K)/\U(p,\K)\times\U(q,\K)$
as the space  $SW_{p,q}=SW_{p,q}(\K)$
of dissipative matrices having block structure (3.2).

Let us represent the matrix (3.2)  in the form
\begin{equation}                  R=
\left(\begin{array}{cc}1&0\\2L&M+N\end{array}\right)
\end{equation}
where
 $M=M^*\gg0, N=-N^*$.  Then
\begin{equation}\frac12(R+R^*)=
\left(\begin{array}{cc}1&L^*\\L&M\end{array}\right)
\end{equation}
Formula  (1.2) implies the identities
\begin{gather}
\det\begin{bmatrix}1&L^*\\L&M\end{bmatrix}_{q-p+j}
=\det[M-LL^*]_j \\
\det\begin{pmatrix}1&L^*\\L&M\end{pmatrix}
=\det(M-LL^*)
\end{gather}

{\sc Lemma 3.1.} {\it Dissipativity of the matrix {\rm (3.3)}
 is equivalent to
positive definiteness of the matrix $M-LL^*$.}

{\sc Proof.} It is enough to apply Sylvester criterion, see (3.5).\kvadrat

\medskip

{\bf 3.2. Minimal parabolic subgroup.}

{\sc Lemma 3.2.} {\it The group of all fractional linear maps of the wedge 
$W_q$ mapping the set $SW_{p,q}$ to itself is isomorphic to the group 
$\U(q-p,\K)\times \U(p,q,\K)$.  Moreover elements of the first factor  
$\U(q-p)$ fix all points of $SW_{p,q}$.}

{\sc Proof.} The statement is equivalent to the following fact.
The group  $H$ of all fractional linear transformations of the matrix ball, 
$\B_{q,q}$ mapping the set of matrices (3.1) to itself, is
$\U(q-p,\K)\times \U(p,q,\K)$. First the inclusion  
 $H\supset\U(q-p,\K)\times \U(p,q,\K)$
is obvious.
Secondly the group   $\U(q-p,\K)\times \U(p,q,\K)$
is a maximal subgroup in $\U(q,q)$. \kvadrat

Let us consider the group $\U(q,q,\K)$ acting on the wedge  $W_q$
and  maximal parabolic subgroup $\cal Q\subset\U(p,q)$ consisting
of transformations
\begin{equation}R\mapsto URU^*+H\end{equation}
where  $U\in\GL(n,\K)$, a $H=-H^*$.
Consider the group
$${\cal S}={\cal Q}\cap (\U(q-p)\times\U(p,q,\K))$$
 which consist of all transformations (3.7)
mapping the set $SW_{p,q}$ to itself. It is easy
show that these transformations have the form
\begin{multline}
\begin{pmatrix} 1&0\\2L&K\end{pmatrix}
                         \mapsto
\begin{pmatrix} A&0\\C&D\end{pmatrix}
\begin{pmatrix} 1&0\\2L&K\end{pmatrix}
\begin{pmatrix} A^*&C^*\\0&D^*\end{pmatrix}
                         +
\begin{pmatrix} 0&-AC^*\\CA^*&Z\end{pmatrix}  =\\
=
\begin{pmatrix} 1&0\\2(DLA^*+CA^*)&DKD^*+CC^*+2DLC^* +Z\end{pmatrix}
\end{multline}
where
$AA^*=1  ;\qquad Z=-Z^*$

By construction, ${\cal S}\subset \U(q-p,\K)\times \U(p,q,\K) $.
But the factor  $\U(q-p,\K)$ fix all points of the set $SW_{p,q}$
Hence we can assume 
${\cal S}\subset  \U(p,q,\K) $.

Formula (3.8) seems cumbersome. Let us give more pleasant description.
The transformation group $\cal S$ is generated by the following transformations

 1. transformations $(L,K)\mapsto (LA^*,K)$ forming the group $\U(q-p,\K)$,

 2. transformations $(L,K)\mapsto (DL, DKD^*)$ forming the group $\GL(p,\K)$,

 3. transformations $(L,K)\mapsto (L,K+Z)$ forming  an Abel group.

 4. transformations $(L,K)\mapsto  (L+C , K+CC^*+2LC^*)$.

The transformations of the type 4 don't form
 group. Nevertheless the transformations of types
 3 and 4 generate a metabelian group of high 2.

It is readily seen, that  $\cal S$ is a maximal parabolic subgroup in
$\U(p,q,\K)$.
We extract a minimal parabolic subgroup $\cal P$ from 
$\cal S$ by the condition: the matrix $D$ is lower-triangular.

Obviously, eigenfunctions of 
the minimal parabolic subgroup $\cal S$ on the space $SW_{p,q}$
have the form
$$\prod_{j=1}^n\det
\begin{bmatrix} 1&L^*\\L&M\end{bmatrix}_{q-p+j}^{\lambda_j-\lambda_{j+1}}
=\det[M-LL^*]_j^{\lambda_j-\lambda_{j+1}}
$$
These functions are restrictions to $SW_q$ of 
  ${\cal Q}$-eigenfunctions on $W_q$.

{\bf 3.3. Functions ${\cal B}_\alpha$} in our coordinates on $SW_{p,q}$ have the form
\begin{equation}
{\cal B}_\alpha(R) =\biggl( \frac{ \det(2(R+R^*))}
                       {|\det(1+R)|^2}\biggr)^\alpha=
  \biggl(\frac{\det(4(M-LL^*)}
               {|\det(1+M+N)|^2}\biggr)^\alpha
\end{equation}
(they also are restrictions to $SW_q$ of functions ${\cal B}_\alpha$ on $W_q$).

\medskip

{\bf 3.4. The spaces $\O(2n+1,\C)/\O(2n+1)$, $\SOS(4n+2)/\U(2n+1)$.}
In Subsection 1.1 these spaces were realized as spaces  
of anticymmetric $(2n+1)\times(2n+1)$ matrices over $\R$, $\C$  with norm  $<1$.
Let us embed the space of  $(2n+1)\times(2n+1)$-matrices to the space of
$(2n+2)\times(2n+2)$-matrices. For this purpose we wright zero column on the left
and zeros row above. Then we apply the construction of Subsection 2.6
and obtain the models described above in Subsection 0.5.

\vspace{22pt}

{\large\bf \S 4. Calculations in Matrix Cones.}

\medskip

\addtocounter{sec}{1}
\setcounter{equation}{0}
%{0}

In these section we  evaluate  integrals (0.1)--(0.2) over symmetric cones
$\GL(n,\R)/\O(n)$, $\GL(n,\C)/\U(n)$, $\GL(n,\H)/\Sp(n)$.
In generally, these 3 series are essentially more simple than the others, see for instance
\cite{FK}). Our calculations in these case also are essentially simplified. 

\medskip 

{\bf 4.1. Separation of variables in the cone of real symmetric matrices.
} Thus we have the integral
\begin{equation}
\int_{T\gg0}\prod_{j=1}^n
\frac{\det [T]_j^{\lambda_j-\lambda_{j+1}}}
{\det(1+[T]_j)^{\sigma_j-\sigma_{j+1}}}
\det T^{-(n+1)/2}dT
\end{equation}
Represent a matrix $T$ as block $(1+(n-1))\times(1+(n-1))$ matrix:
$$T=
\left(\begin{array}{cc}  P&q^t\\q&r \end{array} \right)$$
By  formula (1.2),
\begin{align}
\det T&=\det P\cdot (r-qP^{-1}q^t)\\
\det(1+T)&=\det(1+P) (1+r-q(1+P)^{-1}q^t)
\end{align}
In new coordinates, the domain of integration $T\gg0$ is replaced to
$$P\gg0;\qquad      r-qP^{-1}q^t>0$$
(it is consequence of Sylvester criterion and equality (4.2)).
Using these notations we wright the original integral
in the form
\begin{multline*}
\int\prod_{j=1}^{n-2}
\frac{\det [P]_j^{\lambda_j-\lambda_{j+1}}}
      {\det(1+[P]_j)^{\sigma_j-\sigma_{j+1}}}\cdot
\frac{\det P^{\lambda_{n-1}-(n+1)/2}}
     {\det(1+P)^{\sigma_{n-1}}}\times\\
\times
(r-qP^{-1}q^t)^{\lambda_n-(n+1)/2}\bigl\{1+r-q(1+P)^{-1}q^t\bigr\}^{-\sigma_n}
dq\,dr\,dP
\end{multline*}

Let us change variable
$u=r-qP^{-1}q^t$ (we change  $r$ to $u$, other variables are the same
). We obtain
\begin{multline*}
\int_{P\gg0}\biggl(\prod_{j=1}^{n-2}
\frac{\det [P]_j^{\lambda_j-\lambda_{j+1}}}
      {\det(1+[P]_j)^{\sigma_j-\sigma_{j+1}}}\cdot
\frac{\det P^{\lambda_{n-1}-(n+1)/2}}
     {\det(1+P)^{\sigma_{n-1}}}\times\\
\times\int_{u>0,q\in\R^{n-1}}
u^{\lambda_n-(n+1)/2}\Bigl\{1+u+qP^{-1}q^t-q(1+P)^{-1}q^t\Bigr\}^{-\sigma_n}
dq\,du\biggr)dP
\end{multline*}
Next we reduce the expression in the curly brackets to the form
$$\bigl\{1+u+q(P(1+P))^{-1}q^t\bigr\}$$
Then we change the variable $q=(q_1,\dots, q_{n-1})\in\R^{n-1}$
 to the variable
$$h=(h_1,\dots,h_{n-1})=q(P(1+P))^{-1/2}\in\R^{n-1},$$
 leaving the variables
$P$, $u$ the same. The Jacobian of the substitution
is
\begin{equation}
\det P^{1/2}\det(1+P)^{1/2}.
\end{equation}
Thus integral can be decomposed to a product
\begin{eqnarray}
\int_{P\gg0}\prod_{j=1}^{n-2}
\frac{\det [P]_j^{\lambda_j-\lambda_{j+1}}}
      {\det(1+[P]_j)^{\sigma_j-\sigma_{j+1}}}\cdot
\frac{\det P^{\lambda_{n-1}-n/2}}
     {\det(1+P)^{\sigma_{n-1}-1/2}}dP\times\qquad\qquad\\   \qquad
\times\int_{u>0,h\in\R^{n-1}}
u^{\lambda_n-(n+1)/2}\bigl\{1+u+hh^t\bigr\}^{-\sigma_n}
dh\,du
\end{eqnarray}
The first fuctor (4.5) has a form (4.1) with other parameters 
(the parameter $n$ is replaced to
$n-1$, and also $\sigma$ are changed).
Hence we can apply to integral (4.5) the same arguments.
 It remains to calculate the second factor.

\medskip

{\bf 4.2. Several lemmas.} Now we will give without proof several standard propositions
which will be intensively used below.

\smallskip

{\sc Lemma 4.1.} {                \it
The area of an unit sphere in $\R^k$ is }
\begin{equation}
                S=\frac{2\pi^{k/2}}{ \Gamma(k/2)}
\end{equation}

\smallskip

{\sc Lemma 4.2.} (see. \cite{PBM}, 3.1.2.1)
\begin{equation}
\int\limits_{u>0}
\int\limits_{v>0}
f(u+v)u^{\mu-1}v^{\nu-1}du\,dv=\B(\mu,\nu)
\int\limits_{w>0}
f(w)w^{\mu+\nu-1}dw
\end{equation}

\smallskip

{\sc Lemma 4.3.}
\begin{equation}
\int\limits_{x>0}
\int\limits_{y>0}
(1+x+y)^{-a}x^{b-1}y^{c-1}dx\,dy=
\frac{\Gamma(b)\Gamma(c)\Gamma(a-b-c)}
     {\Gamma(a)}
\end{equation}

(We substitute $z=u, w=u+v$ to the integral (4.8) and integrate in $z$; integral (4.9) can be easily reduced to (4.8)).

{\bf 4.3. Calculation of separated factor.} The factor  (4.6) is
   a partial case of the following integral.

{\sc Lemma 4.4.}
\begin{multline}
\int_{u>0,\,h\in\R^k}
u^{\lambda-k/2-1}(1+u+|h|^2)^{-\sigma} du\,dh=
                      \frac{\pi^{k/2}\Gamma(\lambda-k/2)\Gamma(\sigma-\lambda)}
                {\Gamma(\sigma)}
\end{multline}

{\sc Proof.}
Passing to spherical coordinates (see (4.7)) in the variable
$h=(h_1,\dots,h_{n-1})$, we obtain
\begin{multline}\frac{2\pi^{k/2}}{\Gamma(k/2)}
\int_{u>0}
\int_{\rho>0}
 u^{\lambda-k/2-1}(1+u+\rho^2)^{-\sigma}\rho^{k-1}d\rho\,du=\\
=\frac{\pi^{k/2}}{\Gamma(k/2)}
\int\limits_{u>0}
\int\limits_{v>0}
 u^{\lambda-k/2-1}(1+u+v)^{-\sigma}v^{k/2-1}dv\,du
\end{multline}
The last expression is a partial case of integral (4.9).  \kvadrat

\medskip

{\bf 4.4. Complex and quaternionic cases} are very similar to real case
. Formula  (4.4) for a Jacobian is replaced to
$$\det P^{\dim\K/2} \det (1+P)^{\dim\K/2} $$
Then all is reduced to Lemma 4.4.

\vspace{22pt}

{\large\bf \S 5. Integrals over Wedges of Dissipative Matrices.}

\vspace{22pt}

\addtocounter{sec}{1}
\setcounter{equation}{0}

\medskip

In Subsection 2.3 we realized the spaces $\O(n,n)/\O(n)\times\O(n)$,
$\U(n,n)/\U(n)\times\U(n)$, $\Sp(n,n)/\Sp(n)\times\Sp(n)$
as spaces of dissipative matrices.
Let us represent a matrix  $R$  in the form
$R=T+S$, where $T=T^*$, $S=-S^*$.

\medskip

{\bf 5.1. Evaluation for the spaces $\O(n,n)/\O(n)\times\O(n)$.}
 Consider the integral
\begin{equation}
\int\limits_{T=T^t\gg0\, S=-S^t}
\prod_{j=1}^n
\frac{ \det[T]_j^{\lambda_j-\lambda_{j+1}}    }
     { \det[1+T+S]_j^{\sigma_j-\sigma_{j+1} }  }
    \det T^{-n} dT\,dS
\end{equation}

 Let us represent the matrices $T$, $S$ as block
 $((n-1)+1)\times((n-1)+1)$ matrices:
$$T=
\left(\begin{array}{cc}P&q^t\\q&r \end{array}\right)
                          ;\qquad    S=
\left(\begin{array}{cc}A&-b^t\\b&0 \end{array}\right)
$$
By formula (1.2), we obtain
\begin{gather}\det T= \det P\cdot(r-qP^{-1}q^t);\notag \\
\det(1+T+S)=\det(1+P+A)\cdot(1+r-(q+b)(1+P+A)^{-1}(q^t-b^t))
\end{gather}
Substituting $u= r-qP^{-1}q^t $,
we obtain
\begin{multline}
\int\limits_{P\gg0,A=A^t}   dP\,dA
\biggl(\prod\limits_{j=1}^{n-2}
\frac{  \det[T]_j^{\lambda_j-\lambda_{j+1}}  }
     {   \det[1+P+A]_j^{\sigma_j-\sigma_{j+1}}  }\cdot
\frac{  \det P^{\lambda_{n-1}-n} }
     {\det (1+P+A)^{\sigma_{n-1}}}\times
\\ \! \!\! \!\!
\times \int\limits_{u>0,\, q,b\in\R^{n-1}}  \!\!\!\!\!\!\!\!\!  \!\!\!\!\!
u^{\lambda_n-n}
    \bigl\{1+u+qP^{-1}q^t-(q+b)(1+P+A)^{-1}(q^t-b^t)\bigr\}^{-\sigma_n}
du\,dq\,db\biggr)
\end{multline}
Let us wright the expression in the curly brackets as
\begin{equation}
\biggl\{1+u+
\left(\begin{array}{cc} q&b\end{array}\right)
X
\left(\begin{array}{c} q^t\\b^t\end{array}\right)
\biggr\}\end{equation}
where
\begin{equation}
X=
\left(\begin{array}{cc} P^{-1}-(1+P+A)^{-1} & (1+P+A)^{-1}
\\-(1+P+A)^{-1} & (1+P+A)^{-1}\end{array}\right)
\end{equation}
Unfortunately, the matrix $X$ is not symmetric, hence we wright (5.4)
in the form
$$%\begin{equation*}
\biggl\{1+u+
\frac12\left(\begin{array}{cc} q&b\end{array}\right)
(X + X^t)
\left(\begin{array}{c} q^t\\b^t\end{array}\right)
\biggr\}$$%\end{equation*}

\smallskip

{\sc Lemma 5.1.} a) {\it The matrix $\frac12(X+X^t)$ is positive definite.}

\smallskip

b) $\det (\frac12(X+X^t))= \det P^{-1}\det(1+P+A)^{-2}$

\smallskip

{\sc Proof.}
b)
$\frac12(X+X^t)= $
$$
\left(\begin{array}{cc}
P^{-1}-\frac12 (1+P+A)^{-1} - \frac12 (1+P-A)^{-1} &
\frac12 (1+P+A)^{-1} - \frac12 (1+P-A)^{-1}   \\
\frac12 (1+P-A)^{-1} - \frac12 (1+P+A)^{-1}   &
\frac12 (1+P+A)^{-1} + \frac12 (1+P-A)^{-1}
    \end{array}\right)
$$

Adding the second column to the first row, and the second row to the first row,
we obtain

\begin{eqnarray*} \det(\frac12(X+X^t))=
\left(\begin{array}{cc}
P^{-1}& (1+P+A)^{-1}\\
(1+P-A)^{-1}          &
\frac12 (1+P+A)^{-1} + \frac12 (1+P-A)^{-1}
    \end{array}\right)
\end{eqnarray*}
Applying formula (1.2), we obtain that the original determinant equals to
\begin{multline*}
\det P^{-1}\times\\
\times\det\biggl(
\frac12 (1+P+A)^{-1} + \frac12 (1+P-A)^{-1}
-(1+P-A)^{-1} P (1+P+A)^{-1} \biggr)
=                                                        \\
=\det P^{-1}\det  (1+P+A)^{-1} \det(1+P-A)^{-1}\times\\
\times\det\biggl(   \frac12 (1+P-A) + \frac12 (1+P+A) -P\biggr)\nonumber
\end{multline*}
It remains to notice, that the expression in the big brackets equals 1 and 
 $(1+P+A)^t=1+P-A$. This complete our evaluation.

 a) The expression $\det(1+T+S)$ doesn't vanish.
By formula (5.2),  $u$  and 
(5.4) also don't vanish.  Thus the quadratic summand 
$\left(\begin{array}{cc} q&b\end{array}\right)
X
\left(\begin{array}{c} q^t\\b^t\end{array}\right)$
is nonnegative.  Since the evaluated determinant 
doesn't vanish, we obtain required statement.
\kvadrat

Then we change the variables  $q\in\R^{n-1}$, 
$b\in\R^{n-1}$ to the variable  $h\in\R^{2n-2}$ in  integral (5.3) 
by the formula
\begin{equation}h=
\left(\begin{array}{cc} q&b\end{array}\right)\sqrt{{\textstyle\frac12}(X+X^*)}
\end{equation}
All other variables remain the same.
The Jacobian of the substitution is
$\det P^{1/2}\det(1+P+A)$. This yields that our integral
decomposes to a product of two integrals
\begin{eqnarray}
\int\limits_{P\gg0,A=-A^t}
\prod_{j=1}^{n-2}
\frac{  \det[T]_j^{\lambda_j-\lambda_{j+1}}  }
     {   \det[1+P+A]_j^{\sigma_j-\sigma_{j+1}}  }\cdot
\frac{  \det P^{\lambda_{n-1}-n+1/2} }
     {\det (1+P+A)^{\sigma_{n-1}-1}}dP\,dA\times \\
\times
\int_{u>0,h\in\R^{2n-2}}
u^{\lambda_n-n}(1+u+|h|^2)^{-\sigma_n}du\,dh
\end{eqnarray}
The second factor can be easily evaluated by Lemma 4.4. The first
factor itself has  form (5.1).

\medskip

{\bf 5.2. Separation of variables for the series $\U(n,n)/\U(n)\times\U(n)$.}
Consider integral (0.5) over a  wedge of complex dissipative 
matrices.
Let us represent $T$, $S$ (recall that $T=T^*$, $S=-S^*$)
 as block
$((n-1)+1) \times ((n-1)+1)$ matrices:
$$
T=\left(\begin{array}{cc} P&q^* \\q&r\end{array}\right)
                 ;\qquad  S=
\left(\begin{array}{cc}A&-b^*\\b&ic \end{array}\right)
                                   $$
Repeating considerations of Subsection 5.1, we 
transform our integral to the form
\begin{multline*}
\int\limits_{P\gg0,A=-A^*}\!\!\!\!  \!\!\!\! \!\!\!\!  dP\,dA
                   \Biggl(
\prod_{j=1}^{n-2}
\frac{ \det[P]_j^{\lambda_j-\lambda_{j+1}}  }
   { \det [1\!\!+\!\!P\!\!+\!\!A]_j^{\{\sigma_j-\sigma_{j+1}\|\tau_j-\tau_{j+1}\}}  }
  \cdot
\frac{ \det P^{\lambda_{n-1}-2n}  }
   { \det (1\!\!+\!\!P\!\!+\!\!A)^{\{\sigma_{n-1}\|\tau_{n-1}\}}  }
\times \\
\shoveleft{
\times      \int\limits_{u>0, b\in\C^{n-1}, q\in\C^{n-1}, c\in\R}
u^{\lambda_n-2n}      \times}\\ \times
\Bigl\{1+u+ic+qP^{-1}q^*-(q+b)(1+P+A)^{-1}(q^*-b^*)
\Bigr\}^{\{-\sigma_n\|-\tau_n\}}du\,db\,dq\,dc\Biggr)
\end{multline*}
The expression in curly brackets has the form
$$\biggl\{1+u+ic+
\left(\begin{array}{cc} q&b\end{array}\right)
X
\left(\begin{array}{c} q^*\\b^*\end{array}\right)  \biggr\}
$$
where the matrix  $X$ is the same as above (5.5). Let us change the variable $c$ to
 $s$ by the formula
$$s=c-     \frac{1}{2i}
\left(\begin{array}{cc} q&b\end{array}\right)
       (X^*-X)
\left(\begin{array}{c} q^*\\b^*\end{array}\right)
$$
The Jacobian of the substitution is  1. The expression in the curly brackets
transforms to
$$\biggl\{1+u+is+    \frac{1}{2}
\left(\begin{array}{cc} q&b\end{array}\right)    (X+X^*)
\left(\begin{array}{c} q^*\\b^*\end{array}\right)
\biggr\}
$$
Consider the substitution
$$h=\det({\textstyle\frac12} (X+X^*))^{1/2}$$
Its Jacobian is $\det P |\det(1+P+A)|^2$ (see. Lemma 5.1). Then
 our integral 
decomposes to a product 
of the integrals
\begin{eqnarray}
\int\limits_{P\gg0,A=-A^*}  \!\!
\prod_{j=1}^{n-2}
\frac{ \det[P]_j^{\lambda_j-\lambda_{j+1}}  }
   { \det [1\!\!+\!\!P\!+\!\!A]_j^{\{\sigma_j-\sigma_{j+1}\|\tau_j-\tau_{j+1}\}}  }
  \cdot
\frac{ \det P^{\lambda_{n-1}-2n+1}  }
   { \det (1\!\!+\!\!\!P\!\!\!+\!\!A)^{\{\sigma_{n-1}-1\|\tau_{n-1}-1\}}  } dP\,dA
\times \nonumber \\
\times     \int\limits_{u>0, h\in\C^{2n-2}, s\in\R}
u^{\lambda_n-2n}
\Bigl\{1+u+is+|h|^2
\Bigr\}^{\{-\sigma_n\|-\tau_n\}}du\,dh\,dc         \qquad
\end{eqnarray}

It remains to calculate multiplier (5.9).

{\bf 5.3. Calculation of separated factor.}
Passing to spherical coordinates in the variable
$h\in\C^{2n-2}\simeq\R^{4n-4}$, we transform (5.9) to the form
\begin{eqnarray*}
\frac{2\pi^{2n-2}}{\Gamma(2n-2)}
\int\limits_{u>0, \rho>0,c\in\R}
u^{\lambda_n-2n}(1+u+ic+\rho^2)^{-\sigma_n}
                (1+u-ic+\rho^2)^{-\tau_n}    \rho^{4n-5}d\rho\,du\,dc
=  \\ =
\frac{\pi^{2n-2}}{\Gamma(2n-2)}
\int\limits_{u>0, v>0,c\in\R}
u^{\lambda_n-2n}(1+u+ic+v)^{-\sigma_n}
                (1+u-ic+v)^{-\tau_n}    v^{2n-3}dv\,du\,dc
\end{eqnarray*}
Applying Lemma 4.2 we obtain
\begin{equation}
\frac{\pi^{2n-2}\Gamma(\lambda_n-2n+1)}
     {\Gamma(\lambda_n-1)}
\int\limits_0^\infty dw \int\limits_{-\infty}^\infty dc
\left\{
 (1+w+ic)^{-\sigma_n} (1+w-ic)^{-\tau_n}   w^{\lambda_n-2}
\right\}
\end{equation}
For integration in the variable $c$, we use the formula
 (see. \cite{PBM}, 2.2.5.33) 
$$\int\limits_{-\infty}^\infty
\frac{dx}
  { (a+ix)^\mu(b-ix)^\nu }=
2\pi(a+b)^{1-\mu-\nu}
\frac{\Gamma(\mu+\nu-1)}
  { \Gamma(\mu)   \Gamma(\nu) }$$
(where $a>0$, $b>0$) and obtain
\begin{eqnarray*}
\frac{2\pi^{2n-1}\Gamma(\lambda_n-2n+1)}
     {\Gamma(\lambda_n-1)}
\int\limits_0^{\infty}(2+2w)^{1-\sigma_n-\tau_n} w^{\lambda_n-2}dw
\cdot
\frac{ \Gamma(\sigma_n+\tau_n-1)}
     { \Gamma(\sigma_n)   \Gamma(\tau_n)}
\end{eqnarray*}
Last integral reduces to  $\B$-function in the standard way,
 and finally we obtain
$$
\frac{  2^{2-\sigma_n-\tau_n} \pi^{2n-1}
         \Gamma(\lambda_n-2n+1)\Gamma (\sigma_n+\tau_n-\lambda_n)  }
   { \Gamma (\sigma_n) \Gamma (\tau_n)  }
$$

\smallskip

{\bf 5.4. The case $G/K=\Sp(n,n)/\Sp(n)\times\Sp(n)$.}
Consider integral (0.6) over the wedge of dissipative quaternionic matrices
 $R=T+S$. Recall that  $T=T^*$, $S=-S^*$.
As above, we represent matrices  $T$, $S$
 in the form
$$
T=\left(\begin{array}{cc} P&q^* \\q&r\end{array}\right)
                 ;\qquad  S=
\left(\begin{array}{cc}A&-b^*\\b&c \end{array}\right)
                                   $$
where $c=c_1\i+c_2\j+c_3\k$ is a purely imaginary quaternion.
Repeating considerations of Subsection 5.2, we reduce  integral (0.6) to the form
\begin{eqnarray}
\int\limits_{P\gg0,A=-A^*}
\prod\limits_{j=1}^{n-2}
\frac{  \det[T]_j^{\lambda_j-\lambda_{j+1}}  }
     {   \det[1+P+A]_j^{\sigma_j-\sigma_{j+1}}  }\cdot
\frac{  \det P^{\lambda_{n-1}-4n+2} }
     {\det (1+P+A)^{\sigma_{n-1}-4}}dP\,dA\times \nonumber\\
\times
\int\limits_{u>0,h\in\R^{8n-8},c=-\bar c\in \H}
u^{\lambda_n-4n}\bigl|1+u+s+|h|^2\bigr|^{-\sigma_n}du\,dh\,dc  \qquad
\end{eqnarray}
It remains to evaluate separated factor (5.11).

Passing to spherical coordinates in variable  $h\in\R^{8n-8}$,
 we obtain
\begin{eqnarray*}
\frac{2\pi^{4k-4}}
     { \Gamma(4k-4)}
\int\limits_{u>0,\rho>0,c=-\bar c\in \H}
  u^{\lambda_n-4n} \Bigl|1+u+\rho^2+c\Bigr|^{-\sigma_n}
  \rho^{8n-9}du\,d\rho\,dc= \\               =
\frac{\pi^{4k-4}}
     { \Gamma(4k-4)}
\int\limits_{u>0,v>0,c=-\bar c\in \H}
  u^{\lambda_n-4n} \Bigl|1+u+v+c\Bigr|^{-\sigma_n}
  v^{4n-5}du\,dv\,dc  =\\                   =
\frac{\pi^{4k-4}}
     { \Gamma(4k-4)}
\int\limits_{u>0,v>0,c=-\bar c\in \H}
  u^{\lambda_n-4n} \Bigl((1+u+v)^2+|c|^2\Bigr)^{-\sigma_n/2}
  v^{4n-5}du\,dv\,dc
\end{eqnarray*}

Then we pass to spherical coordinates in variable  $c\in\R^3$,
and get
\begin{eqnarray*}
\frac{4\pi\cdot\pi^{4k-4}}
     { \Gamma(4k-4)}
\int\limits_{u>0,v>0,z>0}
  u^{\lambda_n-4n} \bigl((1+u+v)^2+z^2\bigr)^{-\sigma_n/2}
  v^{4n-5}z^2du\,dv\,dz=\\             =
\frac{2\pi^{4k-3}}
     { \Gamma(4k-4)}
\int\limits_{u>0,v>0,x>0}
  u^{\lambda_n-4n} \bigl((1+u+v)^2+x\bigr)^{-\sigma_n/2}
  v^{4n-5}x^{1/2}du\,dv\,dx
\end{eqnarray*}
Integration in $x$ gives
$$
\frac{2\pi^{4k-3}}
     { \Gamma(4k-4)}
       \B(3/2, \sigma_n/2-3/2)
\int\limits_{u>0,v>0,}
  u^{\lambda_n-4n}v^{4n-5} (1+u+v)^{3-\sigma_n}
  du\,dv $$
Applying Lemma 4.4, we get the final formula
$$
\frac{2\pi^{4k-3}}
     { \Gamma(4k-4)} \cdot
\frac{\Gamma(3/2)\Gamma( \sigma_n/2-3/2 )}
     {\Gamma(\sigma_n/2)  } \cdot
\frac{\Gamma(4k-4) \Gamma(\lambda_n-4n+1)\Gamma(\sigma_n-\lambda_n)}
     {\Gamma(\sigma_n-3)}
$$
Applying Duplication Formula to  $  \Gamma(\sigma_n-3) $,
 we  reduce the last expression to the form
$$
2^{\sigma_n-4}\pi^{4k-3}
\frac{  \Gamma(\lambda_n-4n+1) \Gamma(\sigma_n-\lambda_n)  }
     {\Gamma(\sigma_n/2)  \Gamma(\sigma_n/2-1)}
$$

\medskip

\vspace{22pt}

{\large\bf \S 6. Integrals over the Spaces $\U(p,q,\K)/\U(p,\K)\times\U(q,\K)$.}

\vspace{22pt}

\addtocounter{sec}{1}
\setcounter{equation}{0}

\medskip

In the Section we evaluate integrals (0.7)--(0.9) over the spaces
$\O(p,q)/\O(p)\times\O(q)$, $\U(p,q)/\U(p)\times\U(q)$,
$\Sp(p,q)/\Sp(p)\times\Sp(q)$. For definitenes, we assume
$p\le q$.

{\bf 6.1. The case $G/K=\O(p,q)/\O(p)\times\O(q)$.}
Identities 
(3.5)--(3.6) allow to wright integral (0.7) in the form
\begin{eqnarray}
\int\limits_{M-LL^t\gg0, N=-N^t}
\prod_{j=1}^n
\frac{
\det\left[\begin{array}{cc}1&L^t\\L&M\end{array}\right]_{q-p+j}%
^{\lambda_j-\lambda_{j+1}}    }
     {
\det[1+M+N]_j^{\sigma_j-\sigma_{j+1}}  }
  \det \left(\begin{array}{cc}1&L^t\\L&M\end{array}\right)^{-(p+q)/2}
 dL\,dM\,dN
\end{eqnarray}
Let us represent $M,N$ as  block
 $((p-1)+1) \times((p-1)+1 )$ matrices,  and $L$ as a block 
 $((p-1)+1 ) \times (q-p)$ matrix:
$$
M=\left(\begin{array}{cc}P&q^t\\q&r\end{array}\right);\qquad
N=\left(\begin{array}{cc}A&-b^t\\b&0\end{array}\right)  ;\qquad
L=\left(\begin{array}{c} H\\l\end{array}\right)
$$
By (1.2),
\begin{eqnarray}
\det \left(\begin{array}{cc}1&L^t\\L&M\end{array}\right)=
\det \left(\begin{array}{ccc}
     1&H^t&l^t \\
     H&P&q^t\\
     l&q&r  \end{array}\right) =
\qquad\qquad\qquad\qquad\qquad\qquad\nonumber\\
=\det \left(\begin{array}{cc}1&H^t\\H&P\end{array}\right) \cdot
\biggl\{r -
 \left(\begin{array}{cc}l&q\end{array}\right)
    \left(\begin{array}{cc}1&H^t\\H&P\end{array}\right)^{-1}
     \left(\begin{array}{cc}l^t\\q^t\end{array}\right)
      \biggr\}\nonumber \\
\det(1+M+N)=\det(1\!+\!P\!+\!A)
      \cdot\bigl(1+r-(q+b)(1\!+\!P\!+\!A)^{-1}(q^t-b^t)\bigr) \nonumber
\end{eqnarray}
Denote the expression in curly brackets by $u$.
The integral converts to the form
\begin{eqnarray}
\int dP\,dA\,dH\Biggl(
\prod_{j=1}^{n-2}
\frac{
\det\left[\begin{array}{cc}1&H^t\\H&P\end{array}\right]_{q-p+j}%
^{\lambda_j-\lambda_{j+1}}    }
     {
\det[1+P+A]_j^{\sigma_j-\sigma_{j+1}}  } \cdot
 \frac{ \det \left(\begin{array}{cc}1&H^t\\H&P\end{array}\right)%
^{\lambda_{n-1}-(p+q)/2}}
{\det(1+P+A)^{\sigma_{n-1}}  } \times \nonumber \\
\!\!\!\times\int\limits_{u>0,\,q,b\in\R^{p-1},l\in\R^{q-p}} \!\!\!\!\!\!\!\!\!\!\!
u^{\lambda_p-(p+q)/2}\biggl\{1\!+\!u\! +\!
 \left(\begin{array}{cc}l&q\end{array}\right)
    \left(\begin{array}{cc}1&H^t\\H&P\end{array}\right)^{-1}
     \left(\begin{array}{cc}l^t\\q^t\end{array}\right)+ \\
\!\!\!\!+ \left(\begin{array}{cc}q&b\end{array}\right)
 \left(\begin{array}{cc}
    -(1\!+\!P\!+\!A)^{-1}& -(1\!+\!P\!+\!A)^{-1}\\
      (1\!+\!P\!+\!A)^{-1} & (1\!+\!P\!+\!A)^{-1} \end{array}\right)
     \left(\begin{array}{cc}q^t\\b^t\end{array}\right)
\biggr\}^{-\sigma_p}
du\,dl\,dq\,db\Biggr)
\end{eqnarray}
The expression in the curly brackets in (6.2)--(6.3) has form
\begin{equation}
\biggl\{1+u+
 \left(\begin{array}{ccc}l&q&b\end{array}\right)
                            X
 \left(\begin{array}{c}l^t\\q^t\\b^t\end{array}\right)
\biggr\}
\end{equation}
where
$$X=
     \left(\begin{array}{ccc}
    \left(\begin{array}{cc}1&H^t\\H&P\end{array}\right)^{-1}+
    \left(\begin{array}{cc}0&0\\0&-S^{-1}\end{array}\right)
                        &   \phantom{a}&
 \left(\begin{array}{c}0\\-S^{-1}\end{array}\right)
                \\
  \phantom{d} & & \\
     \left(\begin{array}{cc}0& S^{-1}\end{array}\right)
             &   &
            S^{-1}
\end{array}\right)
$$
and $S$ denote
$$S=1+P+A$$
As in Subsection 5.1, we wright (6.4) in the form
\begin{equation}
\biggl\{1+u+
 \left(\begin{array}{ccc}l&q&b\end{array}\right)
                            \frac12(X+X^t)
 \left(\begin{array}{c}l^t\\q^t\\b^t\end{array}\right)
\biggr\}
\end{equation}

\smallskip

{\sc Lemma 6.1.}
$$\det\left(\frac12(X+X^t)\right)=
\det    \left(\begin{array}{cc}1&H^t\\H&P\end{array}\right)^{-1}
\cdot
\det (1+P+A)^{-2}$$

{\sc Proof.}
$\det\left({\textstyle\frac12}(X+X^t)\right) =$
$$=\det
     \left(\begin{array}{cc}
    \left(\begin{array}{cc}1&H^t\\H&P\end{array}\right)^{-1}+
    \left(\begin{array}{cc}0&0\\0&-\frac12S^{t-1}-\frac12S^{-1}\end{array}\right)
                        &
 \left(\begin{array}{c}0\\-\frac12S^{-1}+\frac12S^{t-1}\end{array}\right)
                \\
  \phantom{d} &  \\
     \left(\begin{array}{cc}0& \frac12S^{-1}-\frac12S^{t-1}\end{array}\right)
             &
           \frac12 S^{-1}+\frac12S^{t-1}
\end{array}\right)
$$
Then we add the third row to the second row and the third column to the second column.
We obtain
$$\det
     \left(\begin{array}{cc}
    \left(\begin{array}{cc}1&H^t\\H&P\end{array}\right)^{-1}
                        &
 \left(\begin{array}{c}0\\(1+P-A)^{-1}\end{array}\right)
                \\
  \phantom{d} &  \\
     \left(\begin{array}{cc}0& (1+P+A)^{-1}\end{array}\right)
             &
           \frac12 (1+P+A)^{-1} +\frac12 (1+P-A)^{-1}
\end{array}\right)
$$
Formula (1.2) reduces the determinant to the form
\begin{eqnarray*}\det
    \left(\begin{array}{cc}1&H^t\\H&P\end{array}\right)^{-1}
      \cdot
\det\Bigl( \frac12 (1+P+A)^{-1}  + \frac12 (1+P-A)^{-1}  - \\-
 \left(\begin{array}{cc}0& (1+P+A)^{-1}\end{array}\right)
  \left(\begin{array}{cc}1&H^t\\H&P\end{array}\right)
 \left(\begin{array}{c}0\\(1+P-A)^{-1}\end{array}\right)
    \Bigr)
=    \\  =
  \det  \left(\begin{array}{cc}1&H^t\\H&P\end{array}\right)^{-1}
\det(1+P+A)^{-1}    \det (1+P-A)^{-1}   \times\\ \times
\det\bigl({\textstyle\frac12} (1+P-A) + {\textstyle\frac12}(1+P+A) -P\bigr)
\qquad\blacksquare
\end{eqnarray*}

After the substitution

$$
 \left(\begin{array}{ccc}l&q&b\end{array}\right)
\sqrt{{\textstyle\frac12}(X+X^t)}=h\in\R^{q-p}\oplus\R^{p-1} \oplus\R^{p-1}
$$
to integral (6.2)--(6.3), we get
\begin{eqnarray*}
\int
\prod_{j=1}^{n-2}
\frac{
\det\left[\begin{array}{cc}1&H^t\\H&P\end{array}\right]_{q-p+j}%
^{\lambda_j-\lambda_{j+1}}    }
     {
\det[1+P+A]_j^{\sigma_j-\sigma_{j+1}}  } \cdot
 \frac{ \det \left(\begin{array}{cc}1&H^t\\H&P\end{array}\right)%
^{\lambda_{n-1}-(p+q)/2+1/2}}
{\det(1+P+A)^{\sigma_{p-1}-1}  }dP\,dA\,dH  \times \nonumber \\
\times\int\limits_{u>0,\,h\in\R^{q+p-2}}
u^{\lambda_p-(p+q)/2}\bigl\{1+u +  |h|^2
\bigr\}^{-\sigma_p}
du\,dh
\end{eqnarray*}
The first factor itself has form (6.1), and the second factor
can be easily evaluated by Lemma 4.4.

\medskip

{\bf 6.2.
The cases $G/K=\U(p,q)/\U(p)\times\U(q)$  and
$\Sp(n,n)/\Sp(n)\times\Sp(n)$.} A separation of variables in these
cases is the same as in the case $G/K=\O(p,q)/\O(p)\times \O(q)$. 
Separated factors coincide with integrals (5.9), (5.11)
evaluated above.

\vspace{22pt}

{\large\bf \S 7. Integrals over the Spaces
  $\Sp(2n,\R)/\U(n)$ 

\hfill and 
 $\Sp(2n,\C)/\Sp(n)$.}

\vspace{22pt}

\addtocounter{sec}{1}
\setcounter{equation}{0}
%{0}

{\bf 7.1. Integrals over $\Sp(2n,\R)/\U(n)$.} We realize this space 
as in Subsection 2.4. Let  $R$ be a complex symmetric matrix,
let $R=T+iS$, where  $T$, $S$ are symmetric. Consider  integral
(0.10).
Let us represent matrices $T$, $S$ in the form
$$T=
  \left(\begin{array}{cc}P&q^t\\q&r\end{array}\right)
, \qquad\qquad S=
  \left(\begin{array}{cc}A&b^t\\b&c\end{array}\right)
                                   $$
Then
\begin{gather*}
\det T=\det P\cdot(r-qP^{-1}q^t)\\
\det (1\!+\!T\!+\!iS)=\det(1\!+\!P\!+\!iA)
\bigl\{1+r+ic -(q+ib) (1+P+iA)^{-1}(q^t+ib^t)\bigr\}
\end{gather*}
Substituting
$$u=r-qP^{-1}q^t$$
we convert (0.10) to the form
\begin{eqnarray*}
\int   dT\,dS \biggl(
\prod_{j=1}^{n-2}
\frac{ \det[P]_j^{\lambda_j-\lambda_{j+1}}}
     {\det[1+P+iA]_j^%
 {\{\sigma_j-\sigma_{j+1}\|\tau_j-\tau_{j+1}\}}   }
 \frac{ \det P^{\lambda_{n-1}-(n+1)}}
  {\det(1+P+iA)^{\{\sigma_{n-1}\|\tau_{n-1}\}} }
\times\nonumber\\
\times
\int u^{\lambda_n-(n+1)}\biggl\{1+u+ic+
  \left(\begin{array}{cc}q&b\end{array}\right)
               X
  \left(\begin{array}{c}q^t\\b^t\end{array}\right)
                             \biggr\}^{\{-\sigma_n\|-\tau_n\}}
du\,dp\,db    \biggr)
\end{eqnarray*}
where
$$X=
  \left(\begin{array}{cc}
     P^{-1}-(1+P+iA)^{-1}&-i(1+P+iA)^{-1}\\
     -i(1+P+iA)^{-1}& (1+P+iA)^{-1}\end{array}\right)
                      $$
Then we repeat the arguments of Subsections 5.1--5.2. 
We consider the substitution
\begin{align*}s=&c+ {\im}
  \left(\begin{array}{cc}q&b\end{array}\right)
               X
  \left(\begin{array}{c}q^t\\b^t\end{array}\right)
               \\
 h=& \left(\begin{array}{cc}q&b\end{array}\right) \sqrt{{\re} X}
\end{align*}
Repeating proof of Lemma 5.1, we obtain
$$\det({\re} X)=
\det P^{-1}\bigl|\det(1+P+iA\bigr|^{-2}$$
Hence the original integral equals  to
\begin{multline}
\int%\limits_{T+S\,\in W(\C)}
\prod_{j=1}^n
\frac{ \det[P]_j^{\lambda_j-\lambda_{j+1}}    }
     { \det[1+P+iA]_j^{\{\sigma_j-\sigma_{j+1}\|\tau_j-\tau_{j+1}\}}  }
\times \\ \times
\frac{ \det[P]_j^{\lambda_{n-1}-(n+1)+1/2}    }
     { \det[1+P+iA]_j^{\{\sigma_{n-1}-1/2\|\tau_{n-1}-1/2\}}  }
dP\,dA\times\\   \times
\int u^{\lambda_n-2n}\bigl\{1+u+is+|h|^2\bigr\}^{\{-\sigma_n\|-\tau_n\}}
du\,ds\,dh
      \end{multline}
Separated factor (7.1) can be easily evaluated by Lemma
4.4.
\medskip

{\bf 7.2. Spaces $\Sp(2n,\C)/\Sp(n)$.}
Consider the realization of this space described in Subsection 2.5.
Represent the matrices $T$, $S$  in the form
$$T=\begin{pmatrix} P&q^*\\q&r\end{pmatrix} \qquad
  S=\begin{pmatrix} A&b^t\\b&c\end{pmatrix}\qquad \mbox{where }
  q,b\in\C^n;\, c\in \C, r\in \R$$
Repeating our calculations, we reduce integral (0.11)
to the form
      \begin{multline}
\int   dP\,dA \biggl(
\prod_{s=1}^{n}
\frac{ \det[P]_s^{\lambda_s-\lambda_{s+1}}}
     {\det[1+P+A\j]_s^%
 {\sigma_s-\sigma_{s+1}}   }
\frac{ \det P^{\lambda_{n-1}-(2n+1)}}
     {\det(1+P+A\j)^%
 {\sigma_{n-1}}   }
\times \\
\times
\int u^{\lambda_n-(2n+1)}\Bigl|1+u+c\j+
\Psi(q,b)\Bigr|^{-\sigma_n}
 du\,dp\,db\,dc    \biggr)
\end{multline}
where
$$\Psi(q,b) =
qP^{-1}q^*-(q+b\j)(1+P+A\j)^{-1}(q^*+\j b^*)$$
Notice that the expression
$${\textstyle\frac12}(\Psi(q,b)-\i^{-1}\Psi(q,b) \i  )$$
is a quaternion having the form $\alpha \j+\beta \k$. Hence we can
change the variable  $c\j=(c_1+c_2\i)\j$ to the variable
$y_1 \j+ y_2 \k$ by the formula
$$y_1 \j+ y_2 \k=c\j+ {\textstyle\frac12}(\Psi(q,b)-\i^{-1}\Psi(q,b) \i  )$$
in integral (7.2).
 Then the expression under symbol of module in  (7.2) changes  to 
\begin{multline}
\Bigl|1+u+  y_1 \j+ y_2 \k +
{\textstyle \frac12}(\Psi(q,b)+\i^{-1}\Psi(q,b) \i  )\Bigr|=\\
= \Bigl|1+u+  y_1 \j+ y_2 \k +qP^{-1}q^* +
\\+{\textstyle\frac12}\Bigl\{
(q+b\j)(1+P+A\j)^{-1}(q^*+\j b^*)  +
(q-b\j)(1+P-A\j)^{-1}(q^*-\j b^*)\Bigr\}\Bigr|
\end{multline}
 The expression in curly brackets has purely complex values
(in spite of the fact that the expression itself contains quaternionic unit $\j$), moreover
this expression is hermitian form depending of complex vector
 $(q, b)$. The determinant of this form 
equals
$$\det P^{-1}\bigl|\det(1+P+A\j)\bigr|^{-2}$$
Thus our integral decomposes to the product
      \begin{multline}
\int
\prod_{s=1}^{n}
\frac{ \det[P]_s^{\lambda_s-\lambda_{s+1}}}
     {\det[1+P+A\j]_s^%
 {\sigma_s-\sigma_{s+1}}   }
\frac{ \det P^{\lambda_{n-1}-(2n-1)+1}}
     {\det(1+P+A\j)^%
 {\sigma_{n-1}-2}   }dP\,dA
\times \\ \times
\int u^{\lambda_n-(2n+1)}
\bigl| 1+u+y_1\j+y_2 \k +|h|^2\bigr|^{-\sigma_n}
du\,dy_1dy_2 dh
\end{multline}
where $h\in\C^{2n-2}$. Factor  (7.3) can be easily evaluated
in the same way as 
 integral (5.11), only integration in $c\in\R^3$
is replaced by integration in  $y\in\R^2$).

\vspace{22pt}

{\large\bf \S 8 Integrals over $\O(2n,\C)/\O(n)$, $\SOS(2n)/\U(n)$.}

\addtocounter{sec}{1}
\setcounter{equation}{0}

\vspace{22pt}

{\bf 8.1. The case  $G/K=\O(2n,\C)/\O(2n)$.}
Denote by $\MR$ the algebra of real $2\times 2$ matrices.
Consider the following elements in $\MR$ 
$$J=\begin{pmatrix}
                0&1\\-1&0\end{pmatrix}
;\qquad\theta= \begin{pmatrix}0&1\\1&0\end{pmatrix}$$
Obviously,
\begin{equation}
J^2=-1,\qquad \theta^2=1;\qquad J\theta =-\theta J
\end{equation}
To any $z\in\C$ we assign
\begin{equation}
 z^\circ=   \alpha+\beta J=
\begin{pmatrix} \alpha&\beta
                \\ -\beta&\alpha
\end{pmatrix}\in \MR
\end{equation}
Matrices $z^\circ$ form a subalgebra
 in $\MR$, which is isomorphic to the algebra
 $\C$.
Let us denote this subalgebra by  $\CC$.
Note that the imaginary unite $i\in\C$ corresponds to the matrix   $J\in\MR$.
 Let us define in $\MR$
also the subalgebra $\RR\subset\CC$ consisting of all matrices of the form $\alpha\cdot 1$,
where $\alpha\in\R$.
It is easy to see that the algebra $\MR$ is generated by the subalgebra $\CC$
and the element $\theta$. For each  $z\in\C$ we have
\begin{equation}
\theta z^\circ = (\overline z)^\circ \theta
\end{equation}

Define  the involution $\Psi\mapsto\Psi^\star$ in $\MR$ by the formula
$$\Psi^\star=
J
\Psi^t
J^{-1}
$$
Equivalently for all $z,w\in\C$
$$ (z^\circ+w^\circ\theta)^\star= {\overline z}^\circ- w^\circ\theta=
\overline{ z}^\circ-\theta {\overline w}^\circ$$
Then
$(pq)^\star=q^\star p^\star$.
We stress that the set of fixed elements of the involution
 $\phantom{}^\star$ coincides with the subalgebra $\R^\circ$.

Let $A$ be an element of the algebra of $k\times k$ matrices over $\MR$.
Let $A_\R$ be the corresponding $2k\times 2k$ matrix over $\R$.  
We define a determinant 
$\DET(A)$ as square root from $|\det(A_{\R})|$.
For matrices over subalgebra $\CC$, this determinant coincides
with module of the usual determinant of a complex matrix.

 Consider the model of the space
  $\O(2n,\C)/\O(2n)$ defined in Subsection 2.6
 (we preserve all notations of  Subsection 2.6)
. Represent a matrix  $R$ in the form
$R=T+H$, where $R=R^t$, $H=-H^t$.
Then we wright $H$  as $H=S\theta$. Thus we represent $R$
in the form
$$R=T+S\theta\qquad
\mbox{where}\quad
T=T^t\gg0,\quad S^t=-\theta S \theta ,\quad JT=TJ,\quad JS=SJ$$

The matrix  $R$ has the size  $2n\times 2n$. We can consider $R$
as  $n\times n$ matrix consisting of 
 $2\times 2$ matrices. In other words, we can consider
$R$ as matrix over $\MR$.    The condition
\begin{equation}
R^t=J^{-1}RJ
\end{equation}
 (see (2.6))
on this language  means that $R$ is $\star$-hermitian
(i.e.  
$R$ as matrix over  $\MR$ coincides with matrix obtained 
by transposition and element-wise involution 
 $\phantom{a}^\star$).

The conditions  $JT=TJ, JS=SJ$
implies that $2\times 2$ blocks of $T$ and $S$ have form (8.1), i.e. 
they are elements of $\CC$.
 Hence we can consider 
$T$ , $S$ as   $n\times n$ matrices over $\CC$. Moreover,
the matrix $T$ is a hermitian matrix over $\CC$,
and $S$   is an antisymmetric matrix over $\CC$.

In these notations integral (0.14) has the form
$$
\int\prod_{j=1}^n
  \frac{ \DET [T]_j^{\lambda_j-\lambda_{j+1}}}
        {\DET [1+T+S\theta]_j^{\sigma_j-\sigma_{j+1}}}
        \DET T^{-(2n-1)} dT\,dS$$
Then we wright $\CC$-matrices $T$, $S$ in the form
$$T=
\begin{pmatrix} P&q^{\circ*}\\q&r\end{pmatrix}
\qquad S=
\begin{pmatrix} A&-b^{\circ t}\\b&0\end{pmatrix},
\qquad\mbox{where}\quad r\in\RR
$$
where symbols  $\phantom{a}^{\circ*}$, $\phantom{a}^{\circ t}$   emphasis
that we consider transposition and conjugation over $\CC$.

Repeating our calculations, we reduce the integral to the form
\begin{multline*}
\int dP\,dA \biggl(
\prod_{j=1}^{n-2}
  \frac{ \DET [P]_j^{\lambda_j-\lambda_{j+1}}}
        {\DET [1+P+A\theta]_j^{\sigma_j-\sigma_{j+1}}}  \cdot
  \frac{ \DET P^{\lambda_{n-1}-(2n-1)}}
        {\DET (1+P+A\theta)^{\sigma_{n-1}}}
                    \times\\ \times
\int u^{\lambda_n-(2n-1)}
\DET(1+u+\Psi(q,b))^{-\sigma_n} du\,dq\,db\biggr)
\end{multline*}
where
 $$
\Psi(q,b)=qP^{-1}q^{\circ*}-(q+b\theta)(1+P+A\theta)^{-1}(q^{\circ*}
-b^{\circ t}\theta))$$
We wright the last expression in two following ways
\begin{align}
\Psi(q,b) & =  qP^{-1}q^{\star}
-(q+b\theta)(1+P+A\theta)^{-1} (q+b\theta)^\star   \\
\Psi(q,b) & = \begin{pmatrix}q&b   \end{pmatrix}
X \begin{pmatrix}q^{\circ*}\\b^{\circ*}   \end{pmatrix}
\end{align}
where
$$X=
\begin{pmatrix}P^{-1}-(1+P+A\theta)^{-1} & -(1+P+A\theta)^{-1}\theta \\
\theta (1+P+A\theta)^{-1} & \theta(1+P+A\theta)^{-1}\theta \end{pmatrix}
$$
It is clear from (8.5), that $\Psi(q,b) \in \RR$
(since $P$ and $(1+P+A\theta)$ are $\star$-hermitian;
this implies that $\Psi(q,b)\in\MR$ is a stable element of involution
 $\star$ ).

\smallskip

{\sc Remark.} The summand  $qP^{-1}q^{\star}$ is not linear over $\MR$.
Hence the expression $\Psi(q,b)$ is not  $\star$-hermitian form 
on module $\MR^{n-1}$ over algebra $\MR$.

\smallskip

We have
$$\Psi(q,b)={\textstyle \frac12}(\Psi(q,b)+ J^{-1}\Psi(q,b) J)$$
Applying this equality to (8.6), we obtain 
$$
\Psi(q,b) = \begin{pmatrix}q&b   \end{pmatrix}
Y \begin{pmatrix}q^{\circ*}\\b^{\circ*}   \end{pmatrix}
 $$
where
$$  Y=
\begin{pmatrix}
 \scriptstyle P^{-1}-\frac12(1+P+A\theta)^{-1}-\frac12(1+P-A\theta)^{-1} &
\scriptstyle -\frac12(1+P+A\theta)^{-1}\theta+\frac12(1+P-A\theta)^{-1}\theta \\
\scriptstyle\frac12\theta (1+P+A\theta)^{-1}-\frac12\theta (1+P-A\theta)^{-1} &
\scriptstyle\frac12 \theta(1+P+A\theta)^{-1}\theta +\frac12 \theta(1+P-A\theta)^{-1}\theta
\end{pmatrix}
$$
First notice that the matrix  $Y$ doesn't depend on  $\theta$,
in spite of expression for $Y$ contains $\theta$. Formally speaking $Y$
 is a matrix over $\CC$. Secondly the matrix $Y$ is
$\star$-hermitian (since $X$, $J^{-1}XJ$ are $\star$-hermitian).
Hence $Y$ is a hermitian matrix over $\CC$.
%
%Sledovatel'no, esli my smenim v formule (8.5)
%simvol $\theta$ na $-\theta$, to vyrazhenie ne izmenitsja. Poetomu
%\begin{multline*}
%  \Psi(q,b)=qP^{-1}q^{\circ*}-\\ -
%{\textstyle \frac12}\bigl\{
%(q+b\theta)(1+P+A\theta)^{-1}(q^{\circ*}-b^{\circ t}\theta)+
%              (q-b\theta)(1+P-A\theta)^{-1}(q^{\circ*}+b^{\circ t}\theta)
%\bigr\}
%\end{multline*}
%No eto vyrazhenie uzhe javljaetsja ermitovoj formoj ot kompleksnyh
%peremennyh $b,q\in\C^{\circ(n-1)}$.
 It is easy to show that its determinant equals
$$\DET(P)^{-1}\DET(1+P+A\theta)^{-2}$$

We obtain
\begin{multline}
\prod_{j=1}^{n-2}
  \frac{ \DET [P]_j^{\lambda_j-\lambda_{j+1}}}
        {\DET [1+P+A\theta]_j^{\sigma_j-\sigma_{j+1}}}
  \frac{ \DET P^{\lambda_{n-1}-(2n-1)+1}}
        {\DET (1+P+A\theta)^{\sigma_{n-1}-2}} \int dP\,dA
                    \times\\ \times
\int\limits_{u>0, h\in\C^{2n-2}} u^{\lambda_n-(2n-1)}
\DET(1+u+|h|^2)^{-\sigma_n}du\,dh
\end{multline}

Factor (8.7) can be easily evaluated by Lemma 4.4.

\medskip

{\bf 8.2. The case $\SOS(2n)/\U(n)$.}
Denote by  $\MC$ the algebra of complex 
  $2\times 2$ matrices over $\C$.
 The algebra $\MC$
contains the subalgebra $\HH$, which consists of matrices having the form
\begin{equation}
\begin{pmatrix}
\alpha&\beta\\
\overline\beta&\overline\alpha
\end{pmatrix}
\end{equation}
Obviously, the subalgebra $\HH$ is isomorphic to the algebra of quaternions$\H$.

For a  
 $k\times k$ matrix $X$ over $\MC$,
we define the determinant $\DET(X)$ as square root from determinant of the same
matrix considered as $2k\times 2k$ matrix over $\C$
(we will use this definition only for dissipative matrices,
hence square root is well-defined).

We preserve notations of Subsection 2.6.  Let us represent the matrix $R$
(having the size $2n\times 2n$ and satisfying the conditions
$R^t=J^{-1}RJ$) as a sum
$R=T+iS$, where $ T=T^*, S=S^*$.
Then
\begin{equation} T=J^{-1}\overline T J\qquad S= J^{-1}\overline  S J
\end{equation}

As in preceding Subsection, we can consider the matrix $R$ as
$n\times n$ matrix over $\MC$. Conditions (8.9)
mean that $T$, $S$ are matrices over  $\HH$.

Let us wright integral (0.15) in the form
$$\int\prod_{j=1}^n
  \frac{ \DET [T]_j^{\lambda_j-\lambda_{j+1}}}
        {\DET [1+T+iS]_j^{\{\sigma_j-\sigma_{j+1}\|\tau_j-\tau_{j+1}\}}}
        \DET T^{-2(2n-1)} dT\,dS$$
where $i$ means multiplication by the scalar $i$ in algebra of complex matrices
having the size $2n\times 2n$. Then we wright the matrices
 $T$, $S$  over $\HH$ in the form
$$T=
\begin{pmatrix} P&q^{\circ*}\\q&r\end{pmatrix}
\qquad S=
\begin{pmatrix} A&b^{\circ*}\\b&c\end{pmatrix},
\quad\mbox{where} \qquad b,q\in\H^{\circ{n-1}}\,\, r,c\in\RR
$$
We get
\begin{multline}\int dP\,dA\biggl(
\prod_{j=1}^{n-2}
  \frac{ \DET [T]_j^{\lambda_j-\lambda_{j+1}}}
        {\DET [1+T+iA]_j^{\{\sigma_j-\sigma_{j+1}\|\tau_j-\tau_{j+1}\}}}
\cdot  \frac{ \DET T^{\lambda_{n-1}-2(2n-1)}}
        {\DET (1+T+iA)^{\{\sigma_{n-1}\|\tau_{n-1}\}}}
   \times\\ \times
\int u^{\lambda_n - 2(2n-1)}\DET(1+u+ic +\Psi(q,b))^{\{-\sigma_n\|-\tau_n\}}
 du\,dc\,db\,dq\biggr)
\end{multline}
where
$$\Psi(b,q)=qP^{-1}q^{\circ*}-(q+ib)(1+P+iA)^{-1}(q^{\circ*}+ib^{\circ*})$$
This expression satisfies to the identity 
$$\Psi(b,q)^t = J^{-1}\Psi(b,q) J$$
where the symbol $\phantom{a}^t$ means transposition 
 in the algebra $\MC$ of $2\times 2$-matrices.
Hence $\Psi(b,q)$ has the form
$\begin{pmatrix}\mu&0\\0&\mu\end{pmatrix}$,
where $\mu\in\C$.

Consider the substitution
$$y=c+ \im (q+ib)(1+P+iA)^{-1}(q^{\circ*}+ib^{\circ*})  $$
(the symbol $\im$ here means an imaginary part of a complex $2\times 2$
matrix). Then  multiplier (8.10) converts to the form
$$\int u^{\lambda_n - 2(2n-1)}\DET(1+u+iy + \re\Psi(q,b))^{\{-\sigma_n\|-\tau_n\}}
     du\,  dc\,db\,dq
         $$
The expression
$\re\Psi(q,b) $ is a hermitian form over  $\HH$,
and its determinant is
$$\DET P^{-1} |\DET(1+P+iA)|^{-2}$$
Hence the integral transforms to the form
\begin{multline}\int dP\,dA\biggl(
\prod_{j=1}^{n-2}
  \frac{ \DET [P]_j^{\lambda_j-\lambda_{j+1}}}
        {\DET [1+P+iA     ]_j^{\{\sigma_j-\sigma_{j+1}\|\tau_j-\tau_{j+1}\}}}
\times
\\
\times  \frac{ \DET P^{\lambda_{n-1}-2(2n-1)+2}}
        {\DET (1+P+iA)^{\{\sigma_{n-1}-2\|\tau_{n-1}-2\}}}
   \times\\ \times
\int\limits_{u>0,h\in\H^{2n-2}}
 u^{\lambda_n - 2(2n-1)}\DET(1+u+ic +|h|^2)^{\{-\sigma_n\|-\tau_n\}}
 du\,dc\,db\,dq\biggr)
\end{multline}
Factor (8.11) can be easily evaluated by Lemma 4.4.

\vspace{22pt}

{\large\bf \S 9. Plancherel Formula for Kernel Representations.}

\nopagebreak

\vspace{22pt}

\nopagebreak

\addtocounter{sec}{1}
\setcounter{equation}{0}

Here we discuss only the case of groups  $G=\O(p,q)$,
considerations for all other series are the same.

\smallskip

For the spaces  $G/K=\O(p,q)/\O(p)\times\O(q)$ we consider 'section of wedge' model
$SW_{p,q}(\R)$, 
see 3.1. We realize a point of the space $G/K$ 
as  a $((q-p)+p)\times ((q-p)+p)$ matrix $R$ having the form
$$R=\begin{pmatrix} 1&0\\2L&M+N\end{pmatrix},\quad \mbox{where}
\quad M=M^t, N=-N^t, M-LL^t\gg 0$$
Let $s_1,\dots,s_p$ be real numbers.
Consider the function $u_s$  on $SW_{p,q}(\R)$
given by the formula
$$u_s(R)=
\prod \det[M-LL^t]_j^{-1/2+i(s_j-s_{j+1})}
$$
Recall that (see 3.2) that these functions are eigenfunctions of
a minimal parabolic subgroup
 ${\cal P}\subset G$. Denote by $\phi_s$ the average of the function
$u_s$ under the action of the maximal compact subgroup 
$K=\O(p)\times\O(q)\subset\O(p,q)$:
$$\phi_s(R)=\int_K u_s(R^{[k]})dk$$
where $dk$ denote Haar measure on  $K$(normalized by the condition:
 the total measure of the group is 1) and $R^{[g]}$ is a fractional
linear transformation defined by the formula (2.3).

It is easy to see that  $\phi_s$ are exactly spherical functions on
the symmetric space $G/K$(see for instance \cite{Hel}, \cite{FK}).

Berezin function ${\cal B}_\alpha$  in our coordinates 
 is given by the formula
$${\cal B}_\alpha(R)=
\frac{4^{\alpha(p+q)}\det(M-LL^t)^{\alpha}}
{\det(1+M+N)^{2\alpha}}$$
(see Subsection 3.3).
Its spherical transform  (see \cite{Hel},\cite{FK})
equals to
$$A_\alpha(s)=\int_{G/K} {\cal B}_\alpha(R) \phi_{-s}(R)dL\,dM\,dN$$
By $K$-invariance of the function ${\cal B}_\alpha(R)$,
this integral equals to the integral
$$\int_{G/K} {\cal B}_\alpha(R)u_{-s}(R)dL\,dM\,dN$$
The last integral is a partial case of integral (0.7).
Thus
$$A_\alpha(s)=\const\cdot\prod_{j=1}^p
  \frac{\Gamma(\alpha+is_j-(q+p)/4+1/2)\Gamma (\alpha-is_j-(q+p)/4+1/2)}
   {\Gamma(2\alpha_j-(m-j))}$$
By  Gindikin--Karpelevich inversion formula (see. \cite{GK},
\cite{Hel})
\begin{equation}
{\cal B}_\alpha(R)= \const
\int_{s_1\ge s_2\ge\dots\ge s_p\ge 0}
 A_\alpha(s)\phi_s(R)\,d\mu(s)
\end{equation}
where $d\mu(s)$ is Gindikin--Karpelevich measure.

Wrighting (9.1) in the explicit form we obtain the following theorem

\medskip

{\sc Theorem 9.1.} {\it Let $\alpha>(q+p)/4-1/2$. Then}
\begin{multline}
{\cal B}_\alpha (R)=
\prod_{k=1}^p\frac{1}{\Gamma(2\alpha-(p-k))}  \times\\ \times
\int_{s_1\ge s_2\ge\dots\ge s_p\ge0}  \quad
\prod_{k=1}^p \biggl|\frac%
  {\Gamma(\alpha+is_k-(q+p)/4+1/2) \Gamma((q-p)/2+is_k)}
  {\Gamma(is_k)}
  \biggr|^2           \times\\ \times
\prod_{1\le k< l\le p}\biggl|
\frac{\Gamma(1/2+i(s_k-s_l))\Gamma(1/2+i(s_k+s_l)) }
     {\Gamma(i(s_k-s_l)) \Gamma(i(s_k+s_l)) }
\biggr|^2\phi_s(R)ds_1\,ds_2\dots ds_p
\end{multline}

\medskip

{\sc Proof.} Formula (9.2) is certainly correct if 
$$\int{\cal B}_\alpha(R)<\infty,\qquad
\int|A_\alpha(s)|\,d\mu(s)<\infty$$
(see \cite{Hel}, Theorem IV.9.5., \cite{FK}, Theorem XIV.5.3).
The both conditions are satisfied for sufficiently large $\alpha$.
By analytical continuation considerations,  formula (9.2)
is correct for all $\alpha>(q+p)/4-1/2$. \kvadrat

{\sc Remark.}
At the point $\alpha=(q+p)/4-1/2$   in the right part of (9.2) there appears a pole in numerator.
 At the same point  ${\cal B}_\alpha$ became not square integrable
on the group $G=\O(p,q)$.
 As it is shown in \cite{Ner4},  for $\alpha<(q+p)/4-1/2$ the spectrum
of a kernel representation is not contained in principle series. By this reason,
Theorem 9.1 for $\alpha<(q+p)/4-1/2$ is not correct.


\begin{thebibliography}{cc}

\bibitem{Ber1}
Berezin F.A. {\it Quantization in complex symmetric domains.}
Izv. Akad. Nauk SSSR, Ser. math., 39, 2, 1362--1402 (1975);
English translation: Math USSR Izv. 9 (1975), No 2, 341--379(1976)

\bibitem{Ber2}
Berezin F.A. {\it On relations between covariant and contravariant
symbols of operators v
for complex classical domains.} Dokl. Akad Nauk SSSR, 241, No 1 (1978);
English translation: Sov. Math. Dokl. 19 (1978), 786--789

\bibitem{vD}
van Dijk H., Hille S.C. {\it Canonical representations related to hyperbolic
spaces},
J.Funct.Anal., 147, 109--139   (1997).

\bibitem{FK}
Faraut J., Koranyi A. {\it Analysis in symmetric cones.}
Oxford Univ.Press, (1994)

\bibitem{Gan}
Gantmaher F.R. {\it Theory of matrices},4-th.ed., Moscow, Nauka(1988);
English translation: Chelsea Publ. Corpor., NY (1959)


\bibitem{Gin}
Gindikin S.G. {\it Analysis on homogeneous spaces}.
Uspehi mat. nauk,19, No 4, 3--92(1964).

\bibitem{Gin2}
Gindikin S.G. {it Invariant distributions in homogeneous domains}
Funkt.Anal i Prilozh. 9 (1975), No 1, 56--58; English translation:
Funct. Anal. Appl. 9 (1975), No.1, 50--52.

\bibitem{GK}
Gindikin S.G., Karpelevich F.I. {\it Plancherel measure for
Riemann symmetric spaces of non-positive curvature.} Dokl. Akad Nauk SSSR,
 145, 252--255(1962)

\bibitem{Hel}
Helgason S. {\it Groups and geometric analysis.}  Acad. Press (1984).

\bibitem{Hua}
Hua Loo Keng, {\it Harmonic analysis of functions of several complex
variables in classical domains
}. Beijing, 1958(Chinese); Russian translation: Moscow, Inostrannaja literatura (1959); English translation: Amer. Math. Soc., Providence (1963)

\bibitem{Ner1}
Neretin Yu.A. {\it Extension of representations of classical groups
to representations of categories 
.} Algebra i analiz, t.3, No 1, 176--202
(1991); English translation in St.Petersburg Math.J.,Vol3(1992),No 1.

\bibitem{Ner2}
Neretin Yu.A.  {\it Categories of symmetries and infinite-dimensional groups}
Oxford University Press (1996);Russion edition: Moscow, URSS(1998).



\bibitem{Ner3}
Neretin Yu.A. {\it Boundary values of holomorphic functions and some
 spectral problems for unitary representations.}
In collection
 {\it Positivity in Lie groups: open problems}, Hilgert J.,
J.D.Lawson, K.-H.Need, Vinberg E.B. (eds.),
Walter de Gruyter, Berlin (1998).

\bibitem{Ner4}
Neretin Yu.A. {\it On separation of spectra in analysis of Berezin kernels}
(to appear)


\bibitem{NO}
Neretin Yu.A., Olshanskii G.I., {\it  Boundary values of holomorphic functions , singular unitary representations of groups 
$\O(p,q)$  and their limits as  $q\to\infty$.}
Zapiski nauchn. semin. POMI RAN 223, 9--91(1995); English translation:
J.Math.Sci. (1997).% Pros'ba ne ubirat' angijskuju ssylku.





\bibitem{OO} Olafsson, G., Orsted, B.,
Bargmann transform for symmetric spaces, (Odense university, preprint,1996).





\bibitem{OZ} Orsted B., Zhang G. {\it Tensor products of analytic continuations
of discrete series}. Can.\ J.\ Math.,49, N 6, 1224--1241(1997)

\bibitem{PBM}
Prudnikov A.P., Brychkov Yu.A., Marichev O.I.
{\it Integrals and series. V.1. Elementary functions.} Moscow, Nauka(1984)
English translation: Gordon and Breach, 1986
\bibitem{Pya}
Piateckij-Shapiro I.I. {\it Geometry of classical domains and
theory of automorphic functions} , Moscow, Fizmatlit (1961);
English translation: {\it Automorphic functions and the geometry of
classical domains}, Gordon and Breach, NY,1969


\bibitem{UU}
Unterberger A., Upmeier H., {\it The Berezin transform and
 invariant differential operators}.
Comm.Matn.Phys.,164, 563--597(1994)

\bibitem{Zha}
Zhang Genkai {\it Berezin transform on line bundles over bounded symmetric domains}, preprint
(to appear in J.Lie theory)

\end{thebibliography}
\end{document}